\documentclass[10pt]{amsart}
\usepackage{times,amsmath,amsbsy,amssymb,amscd,mathrsfs}
\usepackage{booktabs,graphicx,subfigure,epstopdf,wrapfig,chemarrow}

\usepackage{algorithm2e} 
\usepackage{multicol,multirow}
\usepackage{mathtools}
\usepackage[usenames,dvipsnames,svgnames,table]{xcolor}
\usepackage[all]{xy}
\usepackage{wrapfig}
\usepackage{tcolorbox}
\usepackage{enumitem}

\usepackage{tikz,tikz-cd}
\usepackage[utf8]{inputenc}
\usepackage{pgfplots} 
\usepackage{pgfgantt}
\usepackage{pdflscape}
\pgfplotsset{compat=newest} 
\pgfplotsset{plot coordinates/math parser=false}
\newlength\fwidth

\definecolor{myBlue}{rgb}{0.0,0.0,0.55}
\usepackage[pdftex,colorlinks=true,citecolor=myBlue,linkcolor=myBlue]{hyperref}

\usepackage[hyperpageref]{backref}

\usepackage{comment,enumerate,multicol,xspace}

  \newcounter{mnote}
  \let\oldmarginpar\marginpar
    \renewcommand\marginpar[1]{\-\oldmarginpar[\raggedleft\footnotesize #1]%
    {\raggedright\footnotesize #1}}

\newtheorem{theorem}{Theorem}[section]
\newtheorem{lemma}[theorem]{Lemma}
\newtheorem{corollary}[theorem]{Corollary}
\newtheorem{proposition}[theorem]{Proposition}

\newtheorem{example}[theorem]{Example}

\newtheorem{remark}[theorem]{Remark}

\newcommand{\TT}{\texttt{T}}
\newcommand{\dx}{\,{\rm d}x}
\newcommand{\dd}{\,{\rm d}}

\newcommand{\bs}{\boldsymbol}

\DeclareMathOperator*{\sign}{sign}

\DeclareMathOperator*{\spa}{span}

\renewcommand{\div}{\operatorname{div}}
\newcommand{\tr}{\operatorname{tr}}
\newcommand{\alt}{\operatorname{Alt}}

\newcommand{\step}[1]{\noindent\raisebox{1.5pt}[10pt][0pt]{\tiny\framebox{$#1$}}\xspace}

\newcommand{\vertiii}[1]{{\left\vert\kern-0.25ex\left\vert\kern-0.25ex\left\vert #1 
    \right\vert\kern-0.25ex\right\vert\kern-0.25ex\right\vert}}

\newcommand{\Oplus}{\ensuremath{\vcenter{\hbox{\scalebox{1.5}{$\oplus$}}}}}

\newcommand{\prox}{\operatorname{Prox}}

\newcommand{\graysquare}{\tikz \fill[gray!70] (0,0) rectangle (0.18,0.18);}

\begin{document}
\title{Tangential-Normal Decompositions of the Second Family of Finite Element Differential Forms}
\author{Long Chen}%
 \address{Department of Mathematics, University of California at Irvine, Irvine, CA 92697, USA}%
 \email{chenlong@math.uci.edu}%
 \author{Xuehai Huang}%
 \address{Corresponding author. School of Mathematics, Shanghai University of Finance and Economics, Shanghai 200433, China}%
 \email{huang.xuehai@sufe.edu.cn}%

 \thanks{The first author was supported by NSF DMS-2309777 and DMS-2309785. The second author was supported by the National Natural Science Foundation of China Project 12671432.}

\makeatletter
\@namedef{subjclassname@2020}{\textup{2020} Mathematics Subject Classification}
\makeatother

\begin{abstract}
This paper introduces a novel tangential-normal ($t$-$n$) decomposition for the second family of finite element differential forms, presenting a new framework for constructing bases in finite element exterior calculus. The main contribution is the development of a $t$-$n$ basis in which degrees of freedom and shape functions are explicitly dual, a property that streamlines stiffness matrix assembly and enhances the efficiency of interpolation and numerical integration. Additionally, the integration of the well-documented Lagrange element basis supports practical implementation of finite element differential forms in applications.  
\end{abstract}

\maketitle

%\tableofcontents

\section{Introduction}
Finite Element Exterior Calculus (FEEC)~\cite{ArnoldFalkWinthe2006Finite,arnoldFiniteElementExterior2018} is a framework that merges finite element methods with differential geometry and algebraic topology to solve partial differential equations (PDEs), particularly those involving differential forms. It preserves key geometric and topological properties, which ensures stability and accuracy in numerical solutions. In this paper, we provide geometric decompositions of the second family of finite element differential forms on simplices that leverage the traces and the interactions between subsimplices of different dimensions. We then present several bases derived from these geometric decompositions.

In~\cite{ArnoldFalkWinther2009} and the references therein, bases for shape functions and degrees of freedom (DoFs) of finite element differential forms are presented; later work gave systematic basis constructions, scalar--Whitney resolutions, Hodge duality, and unified extension operators~\cite{Licht2022,ChristiansenRapetti2016,BerchenkoKogan2021,Isaac2022}. 
However, these constructions do not in general yield a basis that is directly dual to the particular degrees of freedom used for implementation.
In finite element computation, a basis dual to the DoFs is important for efficient local and global stiffness-matrix assembly. It also simplifies interpolation and numerical integration by directly linking basis functions to DoFs. In~\cite{ChenChenHuangWei2024}, such bases were constructed for edge and face elements. Here, we generalize this construction to the second family of polynomial differential forms $\mathbb{P}_r\Lambda^k(T):=\mathbb{P}_r(T)\otimes\alt^k(\mathscr T^T)$, where $\alt^k(\mathscr T^T)$ is the space of exterior $k$-forms on the tangent space $\mathscr T^T\cong\mathbb{R}^d$, and $\mathbb{P}_r(T)$ is the polynomial space of degree $r\geq1$ on $T$.
%Our results thus establish a foundation for constructing well-conditioned finite element spaces across a range of applications, particularly in fields such as electromagnetics, elasticity, and fluid dynamics. 

%We now give explicit bases of \eqref{eq:decTe} using $t$-$n$ decomposition. 
Fix an orthonormal basis of $\mathbb R^d$. A natural basis of $\alt^k(\mathscr T^T)$ is $\{\dd x_{\sigma}\mid \sigma\in\Sigma(k,d)\}$, where $\Sigma(m,n)$ is the set of increasing maps from $\{1,\ldots,m\}$ to $\{1,\ldots,n\}$. We now give a tangential-normal basis. For an $s$-dimensional subsimplex $e$ of $T$, let $\{\boldsymbol{t}_1^e,\ldots,\boldsymbol{t}_s^e\}$ and $\{\boldsymbol{n}_1^e,\ldots,\boldsymbol{n}_{d-s}^e\}$ be orthonormal bases of the tangent space $\mathscr T^e$ and normal space $\mathscr N^e$, respectively. For $\bs u\in\mathbb R^d$, write $\dd u=\bs u^{\flat}\in\alt^1$. The identity $\alt^k(\mathscr T^e\oplus\mathscr N^e)=\bigoplus_j\alt^{k-j}(\mathscr T^e)\wedge\alt^j(\mathscr N^e)$ gives the following $t$-$n$ basis:
\begin{equation*}%\label{intro:tnalt}
\begin{aligned}
\alt^k(\mathscr T^T)={\rm span}\{\dd t^e_{\sigma}\wedge\dd n^e_{\tau}\mid\, &\sigma\in\Sigma(k-(\ell-s),s),\ \tau\in\Sigma(\ell-s,d-s),\\
&\ell=\max\{s,k\},\ldots,\min\{k+s,d\}\}.
\end{aligned}
\end{equation*}
Take $e(0)$ as the origin and set $e^*:=\{0,\ldots,d\}\setminus e$. For each admissible $\ell$, every basis element is determined by two choices: $k-(\ell-s)$ tangent indices from $e(1:s)$ and $\ell-s$ normal indices from $e^*$. There are $\binom{s}{k-(\ell-s)}$ and $\binom{d-s}{\ell-s}$ choices, respectively. Vandermonde's identity gives
$$
\binom{d}{k}=\sum_{\ell=\max\{s,k\}}^{\min\{k+s,d\}}\underbrace{\binom{s}{k-(\ell-s)}}_{\text{from } e(1:s)}\underbrace{\binom{d-s}{\ell-s}}_{\text{from } e^*}.
$$
The chosen $\ell-s$ normal indices from $e^*$, together with $e$, determine a subsimplex $f=e\sqcup e^*(\tau)\in\Delta_\ell(T)$ with $f\supseteq e$. This gives the following geometric decomposition, proved in Theorem~\ref{thm:Altdec}:
\begin{equation}\label{intro:decTe}
\alt^k(\mathscr T^T)=\Oplus_{\ell=\max\{s,k\}}^{\min\{k+s,d\}}\Oplus_{f\in\Delta_{\ell}(T),\,f\supseteq e}\star_f\alt^{\ell-k}(\mathscr T^e),
\end{equation}
where $\star_f$ is the Hodge star on $f$, and $\star_f\alt^{\ell-k}(\mathscr T^e)\subset\alt^k(\mathscr T^f)\hookrightarrow\alt^k(\mathscr T^T)$.

Using the geometric decomposition of the Lagrange element $$\mathbb P_r(f)=\Oplus_{s=0}^{\ell}\Oplus_{e\in\Delta_s(f)}b_e\mathbb P_{r-(s+1)}(e),$$ where $b_e$ is the scalar bubble polynomial associated with $e$, we obtain the following geometric decomposition of bubble polynomial forms on $f$, where $\ell=\dim f\geq k$:
\begin{equation}\label{intro:Brf}
\mathbb{B}_r\Lambda^k(f)=\Oplus_{s=\ell-k}^{\ell}\Oplus_{e\in\Delta_s(f)}b_e\star_f\mathbb P_{r-(s+1)}\Lambda^{\ell-k}(e).
\end{equation}
This gives the geometric decomposition
\begin{equation}\label{intro:Prkdec}
\begin{aligned}
\mathbb P_r\Lambda^k(T)&=\Oplus_{\ell=k}^d\Oplus_{f\in\Delta_{\ell}(T)}\mathbb B_r\Lambda^k(f)\\
&=\Oplus_{f\in\Delta_k(T)}\mathbb P_r\Lambda^k(f)\oplus\Oplus_{\ell=k+1}^d\Oplus_{f\in\Delta_{\ell}(T)}\mathbb B_r\Lambda^k(f).
\end{aligned}
\end{equation}
Decomposition~\eqref{intro:Brf} corresponds to the horizontal lines, whereas decomposition~\eqref{intro:decTe} corresponds to the vertical lines in Fig.~\ref{fig:tngrid}. Although a geometric decomposition similar to~\eqref{intro:Prkdec} was presented in~\cite{ArnoldFalkWinther2009}, the characterization of the bubble polynomial spaces in~\eqref{intro:Brf} is written differently from that in~\cite{ArnoldFalkWinther2009} and has face indexing and constant exterior-form factors corresponding to the trace-free decomposition in~\cite[(20)]{Isaac2022}. The scalar extension operators differ, so the individual polynomial summands need not coincide; see Remark~\ref{rm:bubblespace}.

Based on these decompositions, we give explicit formulae for the bases and DoFs that define the global finite element space $\mathbb P_r\Lambda^{k}(\mathcal T_h)$ on a conforming triangulation $\mathcal T_h$. We further combine the well-documented Lagrange basis and the new geometric decomposition to establish dual bases for the shape function space and DoFs of the second family of finite element differential forms.

The rest of this paper is organized as follows. Section~\ref{sec:prelim} provides the necessary preliminaries, including simplex and subsimplex structures, barycentric coordinates, Bernstein polynomials, and tangential-normal ($t$-$n$) bases. In Section~\ref{sec:diffform}, we describe extensions of differential forms on simplices and discuss trace operators. Section~\ref{sec:tnbases} presents a geometric decomposition of alternating forms and various $t$-$n$ bases. In Section~\ref{sec:geodecfullpoly}, we present bases for the shape function space and the degrees of freedom associated with a geometric decomposition of the second family of polynomial differential forms.
In Subsection~\ref{sec:geodecBubbleform}, we present geometric decompositions using bubble polynomial differential forms. Finally, we conclude in Section~\ref{sec:conclusion}.

%%%%%%%%%%%%%%%%%%%%
\section{Preliminaries}\label{sec:prelim}
%%%%%%%%%%%%%%%%%%%%

We mainly follow the notation in~\cite{ArnoldFalkWinther2009} with the simplifications presented in~\cite{ChenHuang2024}. We summarize the most important notation and integer indices below:
\begin{itemize}
 \item $\mathbb R^d: d$ is the dimension of the ambient Euclidean space and $d\geq 2$;
\item $\mathbb P_r: r$ is the degree of the polynomial and $r\geq 0$;
%\item  $\alt^k:$ the space of all skew-symmetric $k$-linear forms;
 \item $\Lambda^k: k$ is the degree of the differential form and $0\leq k\leq d$;
\item $\Delta_{\ell}(T): \ell$ is the dimension of a subsimplex $f\in \Delta_{\ell}(T)$ and $0\leq \ell \leq d$.
%\item \RV{$\dx$: the volume form}.
\end{itemize}

\subsection{Simplices and subsimplices}
Let $T \subset \mathbb{R}^{d}$ be a $d$-dimensional simplex with vertices $\texttt{v}_0, \texttt{v}_1, \ldots, \texttt{v}_d$ in general position. We denote by $\Delta(T)$ the set of all subsimplices of $T$, and by $\Delta_{\ell}(T)$ the set of subsimplices of dimension $\ell$, where $0 \leq \ell \leq d$. 

An abstract $d$-simplex $\TT$ is a finite set with $d+1$ elements. The standard (combinatorial) $d$-simplex is $\texttt{S}_d := \{0, 1, \ldots, d\}$. Any abstract simplex $\TT = \{\TT(0), \ldots, \TT(d)\}$ is combinatorially isomorphic to $\texttt{S}_d$ via the map $i \mapsto \TT(i)$. A $d$-simplex $T$ with vertices $\texttt{v}_0, \ldots, \texttt{v}_d$ is a geometric realization of $\TT$ through $\TT(i) \mapsto \texttt{v}_i$, or of $\texttt{S}_d$ via $i \mapsto \texttt{v}_i$. The topological structure is determined by the abstract simplex.

For a subsimplex $f \in \Delta_{\ell}(T)$, we overload the notation $f$ to represent both the geometric simplex and the algebraic subset of $\texttt{S}_d$. Specifically, $f = \{f(0), \ldots, f(\ell)\} \subseteq \texttt{S}_d$ algebraically, while geometrically,
\[
f = {\rm Convex}(\texttt{v}_{f(0)}, \ldots, \texttt{v}_{f(\ell)}) \in \Delta_{\ell}(T)
\]
is the $\ell$-dimensional simplex spanned by the vertices $\texttt{v}_{f(0)}, \ldots, \texttt{v}_{f(\ell)}$. If $f \in \Delta_{\ell}(T)$ for $0 \leq \ell \leq d-1$, then $f^* \in \Delta_{d-\ell-1}(T)$ denotes the subsimplex of $T$ opposite to $f$. Algebraically, treating $f$ as a subset of $\{0, 1, \ldots, d\}$, we have $f^* \subseteq \{0, 1, \ldots, d\}$ such that $f \sqcup f^* = \{0, 1, \ldots, d\}$ with $\sqcup$ being the union operation of two disjoint sets, i.e., $f^*$ is the complement of the set $f$. Geometrically,
\[
f^* = {\rm Convex}(\texttt{v}_{f^*(1)}, \ldots, \texttt{v}_{f^*(d-\ell)}) \in \Delta_{d-\ell-1}(T)
\]
is the $(d-\ell-1)$-dimensional simplex spanned by the vertices not contained in $f$.

Let $F_i$ denote the $(d-1)$-dimensional face opposite vertex $\texttt{v}_i$, i.e., algebraically $F_i = \{i\}^*$. Here, we reserve the capital letter $F$ for a $(d-1)$-dimensional face of $T$. For $f \in \Delta_{\ell}(T)$, we typically use $e \in \Delta_s(f)$, where $s \leq \ell$, to denote a subsimplex of $f$.

An ordering of the vertices determines an orientation of a simplex. We will use $[\cdot ]$ for the ascending ordering. For a subsimplex $f\in \Delta_{\ell}(T)$, we use $[f]$ to denote the ascending order of $f$, i.e., $[f]$ is an ordering of the set $f$ so that $[f](0) < [f](1) < \ldots <[f](\ell)$. 
%\RV{Orientation. Use $[f]$?}

\subsection{Barycentric coordinates and Bernstein polynomials}
For a domain $\Omega \subseteq \mathbb{R}^{d}$ and an integer $r \geq 0$, let $\mathbb{P}_r(\Omega)$ denote the space of real-valued polynomials defined on $\Omega$ of degree less than or equal to $r$. 
% We define $\mathbb{H}_r(\Omega)$ as the space of real-valued homogeneous polynomials of degree $r$ on $\Omega$. 
For simplicity, we write $\mathbb{P}_r = \mathbb{P}_r\left(\mathbb{R}^{d}\right)$. Thus, if the $d$-dimensional domain $\Omega$ has non-empty interior, we have $\operatorname{dim} \mathbb{P}_r(\Omega) = \operatorname{dim} \mathbb{P}_r = \binom{r+d}{d}$. When $\Omega = \{\texttt{v}\}$ is a point, $\mathbb{P}_r(\texttt{v}) = \mathbb{R}$ for all $r \geq 0$, and we set $\mathbb{P}_r(\Omega) = \{0\}$ when $r < 0$. 

For a $d$-dimensional simplex $T$, let $\lambda_0, \lambda_1, \ldots, \lambda_d$ represent the barycentric coordinate functions with respect to $T$. That is, $\lambda_i \in \mathbb{P}_1(T)$ and $\lambda_i(\texttt{v}_j) = \delta_{i,j}$ for $0 \leq i, j \leq d$, where $\delta_{i,j}$ is the Kronecker delta function. The functions $\{\lambda_i\}_{i=0}^d$ form a basis for $\mathbb{P}_1(T)$. Moreover, we have $\sum_{i=0}^d \lambda_i(x) = 1$ and $0 \leq \lambda_i(x) \leq 1$ for all $i = 0, 1, \ldots, d$ and $x \in T$. The subsimplices of $T$ correspond to the zero sets of the barycentric coordinates. Specifically, for $f \in \Delta_{\ell}(T)$, we have
%\[
%f = \{x \in T \mid \lambda_i(x) = 0 \text{ for } i \in f^*\}.
%\]
%In other words,
\begin{equation}\label{eq:lambdaif}
\lambda_i |_f = 0, \quad \forall~i \in f^*.
\end{equation}

We will use multi-index notation $\alpha \in \mathbb{N}^d$, where $\alpha = (\alpha_1, \ldots, \alpha_d)$ with integers $\alpha_i \geq 0$ for $i=1, \ldots, d$. We define $x^\alpha = x_1^{\alpha_1} \cdots x_d^{\alpha_d}$, and $|\alpha| := \sum_{i=1}^d \alpha_i$. We will also use the set $\mathbb{N}^{0:d}$ of multi-indices $\alpha = (\alpha_0, \ldots, \alpha_d)$, where $\lambda^\alpha := \lambda_0^{\alpha_0} \cdots \lambda_d^{\alpha_d}$ for $\alpha \in \mathbb{N}^{0:d}$. For a subsimplex $f \in \Delta_{\ell}(T)$, we treat $f$ as a set of indices and define the bubble polynomial
\[
b_f := \lambda_f := \lambda_{f(0)} \lambda_{f(1)} \ldots \lambda_{f(\ell)} \in \mathbb{P}_{\ell+1}(f).
\]
For a subsimplex $e \in \Delta(T)$, \eqref{eq:lambdaif} implies that $b_f |_e = 0$ if $f\not\subseteq e$, or equivalently, $f\cap e^*\neq \varnothing$. In particular, $b_f |_{\partial f} = 0$, justifying the term ``bubble polynomial''.

We introduce the simplicial lattice~\cite{Chen;Huang:2021Geometric}, whose geometric realization is known as the principal lattice~\cite{nicolaides1972class,nicolaides1973class}. A simplicial lattice of degree $r$ and dimension $d$ is a multi-index set with $d+1$ components and fixed length $r$, defined as
\[
\mathbb{T}_r^d = \left\{ \alpha = (\alpha_0, \alpha_1, \ldots, \alpha_d) \in \mathbb{N}^{d+1} \mid |\alpha|:= \alpha_0 + \alpha_1 + \ldots + \alpha_d = r \right\}.
\]
%An element $\alpha \in \mathbb{T}_r^d$ is called a node of the lattice. 
The Bernstein representation of the space of polynomials of degree $r$ on a simplex $T$ is
\[
\mathbb{P}_r(T) = {\rm span} \left\{ \lambda^\alpha = \lambda_0^{\alpha_0} \lambda_1^{\alpha_1} \ldots \lambda_d^{\alpha_d} \mid \alpha \in \mathbb{T}_r^d \right\}.
\]
In Bernstein form, for $f \in \Delta_{\ell}(T)$, 
\begin{equation*}%\label{eq:PrfBernsteinbasis}
\mathbb{P}_r(f) = {\rm span} \left\{ \lambda_f^\alpha = \lambda_{f(0)}^{\alpha_0} \lambda_{f(1)}^{\alpha_1} \ldots \lambda_{f(\ell)}^{\alpha_\ell} \mid \alpha \in \mathbb{T}_r^\ell \right\}.
\end{equation*}
Through the natural extension defined by the barycentric coordinates, $\mathbb{P}_r(f) \subseteq \mathbb{P}_r(T)$. 

\subsection{Tangential-normal ($t$-$n$) bases}
For a subsimplex $f\in\Delta_\ell(T)$, choose $\ell$ linearly independent tangential vectors $\{\boldsymbol{t}_1^f,\ldots,\boldsymbol{t}_\ell^f\}$ and $d-\ell$ linearly independent normal vectors $\{\boldsymbol{n}_1^f,\ldots,\boldsymbol{n}_{d-\ell}^f\}$. Together, these $d$ vectors form a $t$-$n$ basis of $\mathbb R^d$. Define the tangent and normal spaces of $f$ by
$$\mathscr T^f:=\operatorname{span}\{\boldsymbol{t}_i^f\mid i=1,\ldots,\ell\},\qquad\mathscr N^f:=\operatorname{span}\{\boldsymbol{n}_i^f\mid i=1,\ldots,d-\ell\}.$$
The vectors need not be normalized, and the vectors within either the tangential or the normal family need not be mutually orthogonal. By definition, however $\mathscr T^f\perp\mathscr N^f$.
Within $\mathscr T^f$, let $\{\hat{\boldsymbol{t}}_1^f,\ldots,\hat{\boldsymbol{t}}_\ell^f\}$ be the dual basis, so that $\hat{\boldsymbol{t}}_i^f\in\mathscr T^f$ and $\hat{\boldsymbol{t}}_i^f\cdot\boldsymbol{t}_j^f=\delta_{i,j}$ for $i,j=1,\ldots,\ell$. Similarly, let $\{\hat{\boldsymbol{n}}_1^f,\ldots,\hat{\boldsymbol{n}}_{d-\ell}^f\}$ be the dual basis in $\mathscr N^f$, so that $\hat{\boldsymbol{n}}_i^f\cdot\boldsymbol{n}_j^f=\delta_{i,j}$ for $i,j=1,\ldots,d-\ell$. Since $\mathscr T^f\perp\mathscr N^f$, the combined set $\{\hat{\boldsymbol{t}}_1^f,\ldots,\hat{\boldsymbol{t}}_\ell^f,\hat{\boldsymbol{n}}_1^f,\ldots,\hat{\boldsymbol{n}}_{d-\ell}^f\}$ is dual to $\{\boldsymbol{t}_1^f,\ldots,\boldsymbol{t}_\ell^f,\boldsymbol{n}_1^f,\ldots,\boldsymbol{n}_{d-\ell}^f\}$. If $\delta_{i,j}$ is replaced by $c_i\delta_{i,j}$ with $c_i\neq0$, the two bases are called biorthogonal or scaled dual bases.

Given a subsimplex $f\in\Delta_\ell(T)$ with $0\leq\ell\leq d-1$, we introduce two bases of its normal space. Recall that $F_i$ is the $(d-1)$-dimensional face opposite the $i$-th vertex. Then $f\subseteq F_i$ for $i\in f^*$. Since $\lambda_i|_{F_i}=0$ and $\lambda_i(\texttt{v}_i)=1$, $\nabla\lambda_i$ is orthogonal to $F_i$ and points towards $\texttt{v}_i$.
%Let $\boldsymbol{n}_{F_i}$ be the unit normal vector of $F_i$ with orientation induced from $T$. That is when the signed volume ${\rm vol}(T) > 0$, then $\bs n_{F_i}$ is the outwards normal direction of $\partial T$ and when ${\rm vol}(T) < 0$, it is inwards.
One basis consists of the normal vectors of all $(d-1)$-dimensional faces containing $f$:
\begin{equation*}%\label{eq:normalbasis}
\{\nabla\lambda_i\mid i\in f^*\},
\end{equation*}
which we call the {\it face-normal basis}.
%\RV{The negative sign is chosen so that $-\nabla\lambda_i$ is outward.}
Next, we present a scaled dual basis in $\mathscr N^f$.
%For $f \in \Delta_{\ell}(T)$, where $\ell = 0, 1, \ldots, d-1$, and for $i \in f^*$, let $[f \cup \{i\}]$ represent the $(\ell+1)$-dimensional face in $\Delta_{\ell+1}(T)$ with vertices $\{i, f(0), \dots, f(\ell)\}$ and with orientation given by the ascending ordering.
For a subsimplex $f\in\Delta(T)$, let $\nabla_f$ denote the surface gradient on its tangent plane. The basis
\begin{equation*}%\label{eq:dualnormalbasis}
\{\nabla_{f\cup\{i\}}\lambda_i\mid i\in f^*\}
\end{equation*}
is called the {\it tangential-normal basis}, since $\nabla_{f\cup\{i\}}\lambda_i\in\mathscr T^{f\cup\{i\}}$. Each vector $\nabla_{f\cup\{i\}}\lambda_i$ is orthogonal to $f$ but tangential to the $(\ell+1)$-dimensional face $f\cup\{i\}$.

%\RV{Due to this duality, we also denoted by
%$$
%\widehat{\nabla \lambda}_i := \nabla_{f \cup \{i\}} \lambda_i.
%$$}
%The vector $\nabla \lambda_i \in \mathbb{R}^d$ is parallel to $-\nabla\lambda_i$. Denote its projection onto $\mathscr{T}^f$ as $\nabla_f \lambda_i$. Using this notation, the vector $-\nabla_{f \cup \{i\}} \lambda_i}$ is parallel to the tangential gradient $\nabla_{f \cup \{i\}} \lambda_i$. Therefore, these two bases can also be written as:
%\begin{align*}
%&\text{Rescaled face normal basis: } &\{ -\nabla \lambda_i \mid i \in f^* \}, \\
%&\text{Rescaled tangential-normal basis: } &\{ -\nabla_{f \cup \{i\}} \lambda_i \mid i \in f^* \}.
%\end{align*}

\begin{figure}[htbp]
\subfigure[At vertex $0$: bases $\{-\nabla \lambda_{1}, -\nabla \lambda_{2}, -\nabla \lambda_{3}\}$ and $\{\boldsymbol{t}_{0,1}, \boldsymbol{t}_{0,2}, \boldsymbol{t}_{0,3}\}$.]{
\begin{minipage}[t]{0.425\linewidth}
\centering
\includegraphics[width=5.cm]{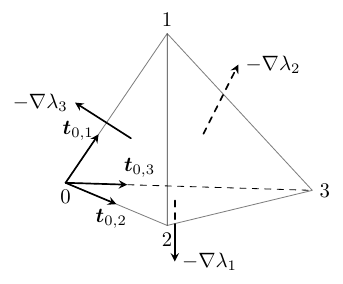}
\end{minipage}}%%
\qquad
\subfigure[On edge $\{0, 1\}$: bases $\{-\nabla \lambda_{2}, -\nabla \lambda_{3}\}$ and $\{\nabla_{[0, 1, 2]} \lambda_2, \nabla_{[0, 1, 3]} \lambda_3\}$.]
{\begin{minipage}[t]{0.435\linewidth}
\centering
\includegraphics*[width=5.3cm]{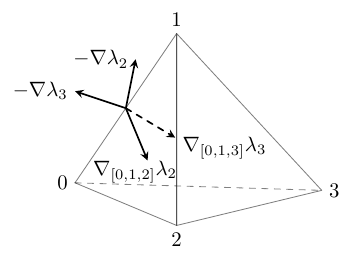}
\end{minipage}}
\caption{Face-normal basis and tangential-normal basis of a vertex and an edge in a tetrahedron. For the face-normal basis functions, a negative sign is added to improve the clarity of the illustration.}
\label{fig:normalbasis}
\end{figure}

%Denote by $\boldsymbol{t}_{i,j}$ the unit edge vector from vertex $\texttt{v}_i$ to vertex $\texttt{v}_j$. By computing the constant directional derivative $\boldsymbol{t}_{i,j} \cdot \nabla \lambda_{\ell}$ using the values at the two vertices, we have
%\[
%\boldsymbol{t}_{i,j} \cdot \nabla \lambda_{\ell} = \delta_{j\ell} - \delta_{i\ell} = 
%\begin{cases} 
%1, & \text{if } \ell = j, \\
%-1, & \text{if } \ell = i, \\
%0, & \text{if } \ell \neq i, j.
%\end{cases}
%\]

%\RV{remove $n^f_{f\cup i}, n_{F_i}$ notation. use $\nabla \lambda$ only. update the figure and examples}.

\begin{example}\label{ex:vertex}\rm
Consider $f \in \Delta_0(T)$, i.e., $f$ is a vertex. Without loss of generality, let $f = \{0\}$. Then $\nabla_{f \cup \{i\}} \lambda_i/|\nabla_{f \cup \{i\}} \lambda_i|$ is the unit tangential vector $\boldsymbol{t}_{0,i}$ from $\texttt{v}_0$ to vertex $\texttt{v}_i$. A scaled dual basis is $\left\{\nabla\lambda_i, i = 1, \ldots, d\right\}$. See Fig.~\ref{fig:normalbasis}~(a).
\end{example}

\begin{example}\rm
Let $f = \{0, 1\}$ be an edge of a tetrahedron. Then we have two bases for the normal plane $\mathscr{N}^f$: $\{\nabla \lambda_{2}, \nabla \lambda_{3}\}$ and $\{\nabla_{[0, 1, 2]} \lambda_2, \nabla_{[0, 1, 3]} \lambda_3\}$. These are scaled dual bases. See Fig.~\ref{fig:normalbasis}~(b).
\end{example}

The following lemma appears in~\cite{ChenHuang2024}. We include it for completeness.
\begin{lemma}\label{lem:normaldual}
For $f\in\Delta_\ell(T)$ with $0\leq\ell\leq d-1$, the tangential-normal basis $$\{\nabla_{f\cup\{i\}}\lambda_i\mid i\in f^*\}$$ of $\mathscr N^f$ is a scaled dual basis of the face-normal basis $\{\nabla\lambda_i\mid i\in f^*\}$.
\end{lemma}
\begin{proof}
Both $\nabla_{f\cup\{i\}}\lambda_i$ and $\nabla\lambda_i$ belong to $\mathscr N^f$ for $i\in f^*$. It suffices to show that $$\nabla_{f\cup\{i\}}\lambda_i\cdot\nabla\lambda_j=0\qquad\text{for }i,j\in f^*,\quad i\neq j.$$ Indeed, $(f\cup\{i\})\subseteq \{j\}^*=F_j$, and hence $\nabla_{f\cup\{i\}}\lambda_i\in\mathscr T^{f\cup\{i\}}\subset\mathscr T^{F_j}\bot \nabla\lambda_j$.
% Since $i\neq j$, this vector is orthogonal to $\nabla\lambda_j$.
\end{proof}

\subsection{Normal bases of an edge on a face}
%Each $\bs n_{e+i,e}$ is the normal vector of $e$ but tangential to $e+i$ for $i\in e^*$. 
%The subspace 
%$$
%\Lambda^{k-\ell}(\mathscr N^e)\wedge \Lambda^{\ell}(\mathscr T^e) \subset \Lambda^{k}(\mathbb R^d).
%$$
%Its element is in the form
%$$
%\dd  n^{\sigma}\wedge \dd  t^{1:\ell}, \quad \sigma \in \Sigma(k-\ell, e^*).
%$$
%Then $\sigma \cup e$ will form a subsimplex $f_{\sigma \cup e}$ of dimension $k$ and $n_{e+i,e}$ is on tangential plane $\mathscr T^f$. So 
%$$
%\Lambda^{k-\ell}(\mathscr N^e)\wedge \Lambda^{\ell}(\mathscr T^e) = \alt^k(\mathscr T^f). 
%$$ 

For an $s$-dimensional subsimplex $e\in\Delta_s(T)$ and an $\ell$-dimensional subsimplex $f\in\Delta_\ell(T)$ with $e\subset f$ and $0\leq s<\ell\leq d$, define the normal space of $e$ in $f$ by $\mathscr N_f^e:=\mathscr N^e\cap\mathscr T^f$.

\vspace{-4pt}
\begin{figure}[htbp]
\centering
\begin{minipage}[t]{0.60\textwidth}
\setlength{\parskip}{0.5em}%
\setlength{\parindent}{0pt}%
\vspace{0pt}

This gives the decomposition
$$
\mathscr T^f=\mathscr T^e\oplus\mathscr N_f^e,
$$
which is a $t$--$n$ decomposition of $\mathscr T^f$ with respect to the lower-dimensional subsimplex $e\in\Delta_s(f)$.

Results for $\mathscr N_T^f=\mathscr N^f$ apply to $\mathscr N_f^e$ after replacing the ambient simplex $T$ with $f$ and the subsimplex $f$ with $e$. Under this replacement, $f^*$ becomes $e^*\cap f=f\setminus e$.
\end{minipage}
\hfill
\begin{minipage}[t]{0.32\textwidth}
\vspace{0pt}
\centering
\includegraphics[width=0.75\linewidth]{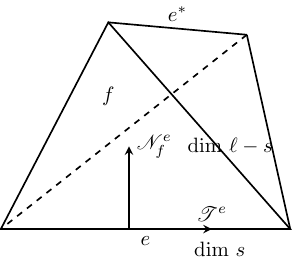}
\vspace{-5pt}
\caption{A $t$--$n$ decomposition of $\mathscr T^f$ based on a subsimplex $e\subset f$.}
\label{fig:Nef}
\end{minipage}
\end{figure}

\begin{lemma}\label{lem:Nfedual}
Let $f\in\Delta_\ell(T)$ and $e\in\Delta_s(f)$ with $0\leq s<\ell\leq d$. Then the two bases
$$
\{\nabla_f\lambda_i\mid i\in f\setminus e\},\qquad\{\nabla_{e\cup\{i\}}\lambda_i\mid i\in f\setminus e\}
$$
of $\mathscr N_f^e$ are scaled dual bases.
\end{lemma}
Moreover, $\{\nabla\lambda_i\mid i\in f\setminus e\}$ is also scaled dual to $\{\nabla_{e\cup\{i\}}\lambda_i\mid i\in f\setminus e\}$, but $\nabla\lambda_i$ need not belong to $\mathscr T^f$. The surface gradient $\nabla_f\lambda_i$ is the projection of $\nabla\lambda_i$ to $\mathscr T^f$.

\subsection{Background on differential forms}
This subsection collects the notation and basic operations for differential forms used throughout the paper. We define the spaces $H\Lambda^k(\Omega)$ in which the finite element spaces will be constructed.

\subsubsection{Increasing sequences}
%For non-negative integers $a$, $b$, $l$, and $m$, with $0\leq b-a\leq m-l$,
For two non-negative integers $m\leq n$, define
$$
\Sigma(m,n):=\{\sigma:\{1,\ldots,m\}\to\{1,\ldots,n\}\mid\sigma(1)<\sigma(2)<\cdots<\sigma(m)\}.
$$
We use $\sigma$ to denote both the map and its range. Thus, for $\sigma\in\Sigma(m,n)$, $\sigma$ also denotes the ordered set $\{\sigma(i)\mid i=1,\ldots,m\}$. The complementary map $\sigma^c\in\Sigma(n-m,n)$ is defined by
$$
\sigma\sqcup\sigma^c=\{1,\ldots,n\}.
$$
%The set \(\Sigma(0: k, 0: d)\) is primarily used to describe subsimplices. For \(\sigma \in \Sigma(0: k, 0: d)\), the subsimplex \(f_\sigma \in \Delta_k(T)\) is formed by the vertices with indices \(\{\sigma(0), \ldots, \sigma(k)\}\). Conversely, for \(f \in \Delta_k(T)\), the indices of its vertices can be sorted in the ascending order to produce an increasing sequence \(\sigma = [f] \in \Sigma(0: k, 0: d)\).
%
%For \(\sigma \in \Sigma(0: k, 0: d)\), the complementary map \(\sigma^* \in \Sigma(1: d-k, 0: d)\) is characterized by
%\[
%\sigma \cup \sigma^* = \{0, 1, \ldots, d\}.
%\]

%For differential forms in \(\mathbb{R}^d\), we use the set \(\Sigma(1: k, 1: d)\), and will abbreviate it as \(\Sigma(k, d) := \Sigma(1: k, 1: d)\).

\subsubsection{Differential forms}
Let $\Omega\subset\mathbb R^d$ be a $d$-dimensional domain. Choose Cartesian coordinates and write $x=(x_1,\ldots,x_d)\in\Omega$. We identify a point $x$ with the vector $\bs x=\vec{ox}\in\mathbb R^d$. Let $\partial_{x_i}$ denote the $i$-th coordinate vector, viewed as an element of $T_o\Omega$. The corresponding dual basis in $(\mathbb R^d)^*$ is $\{\dd x_i\}_{i=1}^d$, where $\dd x_i(\partial_{x_j})=\delta_{i,j}$ for $i,j=1,\ldots,d$. The standard inner product on $\mathbb R^d$ induces an inner product on $(\mathbb R^d)^*$ such that $\langle\dd x_i,\dd x_j\rangle=\delta_{i,j}$ for $i,j=1,\ldots,d$. We reserve $\{\dd x_i\}_{i=1}^d$ for the positively oriented orthonormal basis induced by the ambient coordinate system of $\mathbb R^d$. A general basis is denoted by $\{\dd y_i\}_{i=1}^d$ and need not be orthonormal.

%We can find another basis $\{\dd \hat{y}*i\\}*{i=1}^d$ dual to $\{\dd y_i\}*{i=1}^d$, such that $\\langle \\dd \hat{y}i, \\dd y\_j\\rangle =  \\delta{i,j}$. Specifically, let $M =( \\langle \\dd y\_i, \\dd y\_j \\rangle)*{i,j=1}^d$. Then,
%\[
%(\dd \hat{y}_1, \ldots, \dd \hat{y}_d)^{\intercal} = M^{-1}(\dd y_1, \ldots, \dd y_d)^{\intercal}.
%\]
%When $\{\dd y_i\}_{i=1}^d$ is orthonormal, we have $\dd\hat{y}_i = \dd y_i$ for $i=1,\ldots, d$, as $M$ is the identity matrix.

For a vector space $V$, let $\alt^k(V)$ denote the space of alternating $k$-linear functionals on $V^k:=\underbrace{V\times\cdots\times V}_{k}$. We write $\alt^k$ when $V$ is clear from the context.

For $\omega\in\alt^p$ and $\eta\in\alt^q$, define the {\rm wedge product} $\omega\wedge\eta\in\alt^{p+q}$ by
\begin{eqnarray*}
(\omega\wedge\eta)(v_1,\ldots,v_{p+q})=\sum_{\sigma}{\rm sign}(\sigma)\omega(v_{\sigma(1)},\ldots,v_{\sigma(p)})\eta(v_{\sigma(p+1)},\ldots,v_{\sigma(p+q)}),
\end{eqnarray*}
where the sum is over all permutations $\sigma$ of $\{1,\ldots,p+q\}$ satisfying $\sigma(1)<\cdots<\sigma(p)$ and $\sigma(p+1)<\cdots<\sigma(p+q)$, and ${\rm sign}(\sigma)$ is the signature of $\sigma$.

Let $T_x\Omega$ be the tangent space at $x$. For a smooth manifold $\Omega$, a differential $k$-form is a smooth section of the bundle $\bigcup_{x\in\Omega}\alt^k(T_x\Omega)$. The space of all differential $k$-forms is denoted by $\Lambda^k(\Omega)$, or simply $\Lambda^k$. Since $\Omega$ is a domain in $\mathbb R^d$, translation canonically identifies $T_x\Omega$ with $T_o\Omega$ for every $x\in\Omega$. Thus, the basis $\{\dd y_i\}$ can be used for $\alt^1(T_x\Omega)$ at every $x\in\Omega$.

For $\sigma\in\Sigma(k,d)$, extend the multi-index notation by
$$
\dd y_\sigma:=\dd y_{\sigma(1)}\wedge\cdots\wedge\dd y_{\sigma(k)}.
$$
% Notice that $\{\dd \hat{y}*{\\sigma}, \\sigma\\in \\Sigma(k, d)\\}$ is dual to $\\{\\dd y*{\sigma}, \\sigma\\in \\Sigma(k, d)\\}$, i.e., $\langle\dd \hat{y}*{\\sigma}, \\dd y*{\sigma^c}\rangle=\delta_{\sigma,\sigma^c}$ for $\sigma, \sigma^c \in \Sigma(k, d)$.
An element $\omega\in\Lambda^k(\Omega)$ has the representation
\begin{equation}\label{eq:representation}
\omega=\sum_{\sigma\in\Sigma(k,d)}a_\sigma(x)\dd y_\sigma,\qquad x\in\Omega.
\end{equation}

Using~\eqref{eq:representation}, define the exterior derivative $\dd=\dd_\Omega:\Lambda^k(\Omega)\rightarrow\Lambda^{k+1}(\Omega)$. For $\omega=\sum_{\sigma\in\Sigma(k,d)}a_\sigma(x)\dd y_\sigma$, set
$$
\dd\omega=\sum_{\sigma\in\Sigma(k,d)}\sum_{i\notin\sigma}\partial_{y_i}a_\sigma(x)\dd y_i\wedge\dd y_\sigma.
$$
This definition is independent of the choice of basis. For a subsimplex $f$, denote the exterior derivative on $f$ by $\dd_f:\Lambda^k(f)\to\Lambda^{k+1}(f)$.

For $\omega\in\alt^k(V)$ and $v\in V$, define the {\rm contraction} $\omega\, \lrcorner\, v\in\alt^{k-1}(V)$ by
$$
(\omega\, \lrcorner\, v)(v_1,\ldots,v_{k-1}):=\omega(v,v_1,\ldots,v_{k-1}).
$$
In coordinates, if $\omega=\sum_{\sigma\in\Sigma(k,d)}a_\sigma\dd y_\sigma\in\Lambda^k(\Omega)$, then
$$
\omega\, \lrcorner\, v=\sum_{\sigma\in\Sigma(k,d)}a_\sigma\sum_{i=1}^k(-1)^{i-1}\dd y_{\sigma(i)}(v)\dd y_{\sigma(1)}\wedge\cdots\wedge\widehat{\dd y_{\sigma(i)}}\wedge\cdots\wedge\dd y_{\sigma(k)},
$$
where $\widehat{\dd y_{\sigma(i)}}$ means that $\dd y_{\sigma(i)}$ is omitted.

Since $\sum_{i=0}^d\lambda_i=1$, we have $\sum_{i=0}^d\dd\lambda_i=0$. Thus, $\{\dd\lambda_i\mid i=0,\ldots,d\}$ is linearly dependent. Omitting $\dd\lambda_0$ gives the basis $\{\dd\lambda_i\mid i=1,\ldots,d\}$ of $\alt^1(\mathscr T^T)$. More generally, for $f\in\Delta_\ell(T)$, set
$$
\dd\lambda_f:=\dd\lambda_{f(0)}\wedge\cdots\wedge\dd\lambda_{f(\ell)}\in\alt^{\ell+1}(\mathscr T^T).
$$
Define $\Delta_\ell^0(T) = \{f\in\Delta_\ell(T)\mid 0\in f\}$. For these faces, we use increasing vertex order, so $f(0)=0$ and $[f\setminus\{0\}] = f(1:\ell)$. Then
\begin{equation*}%\label{eq:Altbasis}
\{\dd\lambda_{[f\setminus\{0\}]}\mid f\in\Delta_\ell^0(T)\}\text{ is a basis of }\alt^\ell(\mathscr T^T),
\end{equation*}
where $[f\setminus\{0\}]$ is ordered increasingly. By the wedge-product formula, Lemma~\ref{lem:normaldual}, and Example~\ref{ex:vertex}, a scaled dual basis is
$$
\{\mathop{\textstyle\bigwedge}\nolimits_{j\in [f\setminus\{0\}]}\dd_{[0,j]}\lambda_j\mid f\in\Delta_\ell^0(T)\}.
$$

\subsubsection{Vector proxies}
For a $1$-form $\omega=\sum_{i=1}^d u_i\dd x_i\in\Lambda^1$, define
$$
\prox_1(\omega):=\omega^\sharp=\bs u=(u_1,u_2,\ldots,u_d)^{\intercal},
$$
where $\sharp:\Lambda^1\to\mathbb R^d$ is the sharp operator~\cite{Lee2013}.
%For a vector $\bs t = (t_1, t_2,\ldots, t_n)^{\intercal}$ representing the tangent vector $ \sum_{i=1}^n t_i\partial_{x_i}$, the action is
%$$
%\omega (\bs t) = \sum_{i,j=1}^d u_it_j \langle \dd x_i, \partial_{x_j} \rangle = \bs u\cdot \bs t.
%$$
The map $\prox_1$ is a bijection: if $\bs u=(u_1,\ldots,u_d)^{\intercal}$ and $\omega=\sum_{i=1}^d u_i\dd x_i\in\Lambda^1$, then $\prox_1(\omega)=\bs u$. To match differential-form notation, define
\begin{align*}
\dd u:=\bs u^\flat=\sum_{i=1}^d u_i\dd x_i,
\end{align*}
where $\flat:\mathbb R^d\to\Lambda^1$ is the flat operator~\cite{Lee2013}, and $u$ is the linear $0$-form induced by the inner product, namely $u(\bs x):=(\bs u,\bs x)=\sum_{i=1}^d u_i x_i$ for $\bs x\in\mathbb R^d$.
%When $\bs v\in \mathscr T^f$, $\dd v\in \alt^1(\mathscr T^f)$.
% Here in $\dd u$, the symbol $\dflat$ is understood as an operator that maps a tangent vector $\bs u$ to a cotangent vector $\dd u \in \Lambda^1$. The symbol $\dflat$ is not associated with differentiation. A standard notation for this operation is simply $\flat$, but this lacks a direct connection to differential forms.
% For example, the expression $\dd t \wedge \dd n$ is more naturally interpreted as a 2-form than $\bs t^{\flat} \wedge \bs n^{\flat}$.

By definition, $\dd\lambda_i=(\nabla\lambda_i)^\flat\in\alt^1(\mathscr T^T)$ and $\dd_f\lambda_i=(\nabla_f\lambda_i)^\flat\in\alt^1(\mathscr T^f)$, where $\nabla_f$ is the surface gradient.
%For a given $f\in \Delta_{\ell}(T)$, we introduce notation for the dual basis
%$$
%\dd \hat{\lambda}_i := \dflat (\widehat{\nabla \lambda}*i) = (\\nabla*{f\cup\{i\}} \lambda_i)^{\flat}, \quad i\in f^*.
%$$
%Then $\{\dd \lambda_i \mid i\in f^*\\}$ and $\\{\\dd \\hat{\\lambda}\_i \\mid i\\in f^*\}$ are scaled dual bases of $\alt^{1}(\mathscr N^f)$. Similarly a $f\in \Delta_{\ell}(T)$ and $e\in \Delta_s(f)$
%$$
%\dd_f \hat{\lambda}_i := \dflat (\widehat{\nabla_f \lambda}*i) = (\\nabla*{e\cup\{i\}} \lambda_i)^{\flat}, \quad i\in e^*\cap f = f\setminus  e.
%$$
%And $\{\dd_f \lambda_i \mid i\in f\setminus  e\}$ and $\{\dd_f \hat{\lambda}_i \mid i\in f\setminus  e\}$ are scaled dual bases of $\alt^{1}(\mathscr N^e_f)$. Indeed $\dd_f \hat{\lambda}_i = \dd \hat{\lambda}_i$ for $i\in e^*\cap f = f\setminus  e$.}
%
%\RV{The notation $\dd \hat{\lambda}_i$ in the identity $\dd_f \hat{\lambda}_i = \dd \hat{\lambda}_i$ is potentially ambiguous, as it is unclear whether the differential is taken with respect to $f$ or, for instance, $e$.}
%
By the wedge-product formula and Lemma~\ref{lem:Nfedual}, the following volume forms are scaled dual bases of $\alt^{\ell-s}(\mathscr N_f^e)$:
\begin{equation}\label{eq:dflambdaf-e}
\dd_f\lambda_{[f\setminus e]}:=\mathop{\textstyle\bigwedge}\nolimits_{j\in[f\setminus e]}\dd_f\lambda_j,
\end{equation}
and
\begin{equation}\label{eq:dualdflambdaf-e}
\dd_f\widehat{\lambda}_{[f\setminus e]}:=\mathop{\textstyle\bigwedge}\nolimits_{j\in[f\setminus e]}\dd_{e\cup\{j\}}\lambda_j.
\end{equation}
In both products, the wedge factors are ordered increasingly in $j$.
The pairing of these two elements is
$$
\begin{aligned}
\left\langle\dd_f\lambda_{[f\setminus e]},\dd_f\widehat{\lambda}_{[f\setminus e]}\right\rangle
%&=\det\left(\left\langle\dd_f\lambda_{[f\setminus e](i)},\dd_{e\cup\{[f\setminus e](j)\}}\lambda_{[f\setminus e](j)}\right\rangle\right)_{i,j=1}^{\ell-s}\\
&=\prod_{j\in f\setminus e}\left\langle\dd_f\lambda_j,\dd_{e\cup\{j\}}\lambda_j\right\rangle  =\prod_{j\in f\setminus e}\left|\dd_{e\cup\{j\}}\lambda_j\right|^2.
\end{aligned}
$$
We abbreviate $\dd_f\widehat{\lambda}_{[f\setminus e]}$ by $\dd\widehat{\lambda}_{[e^*]}$ when $f=T$.

%\RV{What is the good notation for its dual basis?}
\subsubsection{Hodge star}
The Hodge star for a positively oriented orthonormal basis $\{\dd x_i \}_{i=1}^d$ is a {\rm linear} operator $\star\,: \alt^k\mapsto \alt^{d-k}$ defined by
\begin{equation}\label{eq:Hodge}
\star(\dd x_{\sigma})
=
{\rm sign} (\sigma,\sigma^c) \dd x_{\sigma^c}, 
\quad
\forall~\sigma \in \Sigma(k, d).
\end{equation}
Here ${\rm sign} (\sigma,\sigma^c)$ is the signature of the permutation $(\sigma(1), \ldots, \sigma(k), \sigma^c(1), \ldots, \sigma^c(d-k))$.
By definition, we have
\[
\star (\star \, \omega) = (-1)^{k(d-k)} \omega, \quad \forall~\omega \in \alt^k, 0\leq k\leq d.
\]
It can be shown that the Hodge star operator $\star$ is an isometry \cite[p.~12]{ArnoldFalkWinthe2006Finite}.

\subsubsection{Inner product}
An intrinsic definition of an inner product on $\alt^k$ is given by
\begin{equation*}
\langle \omega, \eta \rangle := \sum_{\sigma \in \Sigma(k, d)}\omega(e_{\sigma(1)}, \ldots, e_{\sigma(k)})\eta(e_{\sigma(1)}, \ldots, e_{\sigma(k)}),
\end{equation*}
where $(e_{1}, \ldots, e_{d})$ is any orthonormal basis of the tangent space.

Recall that $\{\dd x_i\}_{i=1}^d$ is an orthonormal basis of $\alt^1$, meaning $\langle \dd x_i, \dd x_j \rangle = \delta_{i,j}$ for $i,j = 1,\ldots,d$. This naturally extends to an orthonormal basis $\{\dd x_{\sigma}, \sigma \in \Sigma(k, d)\}$ for $\alt^k$, such that
\begin{equation*}%\label{eq:dualkform}
\langle \dd x_{\sigma}, \dd x_{\tau} \rangle  = \delta_{\sigma, \tau}, \quad \sigma, \tau \in \Sigma(k, d).
\end{equation*}
Denote the volume form by $\dx = \dx_1\wedge \dx_2 \wedge \ldots \wedge \dx_d$, which is positive.
By definition, we have
\begin{equation}\label{eq:innerproductHodge}
\omega \wedge \star \eta = \langle \omega, \eta \rangle \, \dd x, \quad \omega, \eta \in \alt^k.
\end{equation}

%Then
%\begin{equation*}
%\star(\omega \wedge \star \eta) = \langle \omega, \eta \rangle, \quad \omega, \eta \in \alt^k.  
%\end{equation*}

For $\omega, \eta \in \Lambda^k(\Omega)$, a further integral over the domain is included, giving
\[
(\omega, \eta)_{\Omega} = \int_{\Omega}\langle \omega, \eta \rangle \, \dd x = \int_{\Omega} \omega \wedge \star \eta, \quad \omega, \eta \in \Lambda^k(\Omega) .
\]
%For a submanifold $f$ of $\Omega$, the volume form $\dd x$ induces a volume form on $f$, denoted by $\dd x_f$. We define the inner product on $f$ as
%$$
%(\omega, \eta)_f = \int_f \langle \omega, \eta \rangle \, \dd x_f, \quad \omega, \eta \in \Lambda^k(f)\hookrightarrow \Lambda^k(\Omega).
%$$
%\RV{Check if it is consistent with the one defined by pull back.} \RV{Not clear, as the degrees of the volumes $\dd x$ and $\dd x_f$ are different.} \RV{Then define the integral by the pull back of map $f\to T$.}
We denote by $L^2\Lambda^k(\Omega)$ the space of differential forms with square-integrable coefficient functions. The Sobolev space $H\Lambda^k(\Omega)$ is defined as
\[
H\Lambda^k(\Omega) := \{\omega \in L^2\Lambda^k(\Omega) : \dd \omega \in L^2\Lambda^{k+1}(\Omega)\}.
\]

\section{Differential Forms on Simplices}\label{sec:diffform}
In this section, we focus on differential forms on simplices and their relationships across subsimplices. We denote by $\mathbb{P}_r\Lambda^k(T):=\mathbb P_r(T)\otimes \alt^k(\mathscr T^T)$ the space of differential forms on $T$ with polynomial coefficients of degree $r$. 
%The dimension identity
%$$
%\dim \notag \mathbb P_r\Lambda^{k}(T) = {d + r \choose d} {d \choose k} = {d + r \choose r + k} {r +k  \choose k} = { d+r \choose d-k+r} {d-k +r  \choose r}
%$$ 
%can be proved by looking at the coefficients of  \( x_1^{k} x_2^{d-k} x_3^r \) in the expansion of \((x_1 + x_2 + x_3)^{d+r}\) in different ways (group $x_1+x_2$, $x_1+x_3$ or $x_2+x_3$).

\subsection{Extension}
Let $f \in \Delta_{\ell}(T)$ with \(0 \leq k \leq \ell \leq d\). 
Let \(\Pi_f \) be the orthogonal projection \(\mathscr{T}^T = \mathbb R^d \to \mathscr{T}^f\). Then the pullback $\Pi_f^*$ defines an embedding \(\alt^k(\mathscr{T}^f) \hookrightarrow \alt^k(\mathscr{T}^T)\): for $\omega \in \alt^k(\mathscr{T}^f)$, $\Pi_f^*\omega\in \alt^k(\mathscr{T}^T)$ satisfies
\[
\Pi_f^*\omega (v_1, v_2, \ldots, v_k) =\omega(\Pi_f v_1, \Pi_f v_2, \ldots, \Pi_f v_k), \quad v_i \in \mathbb{R}^d=\mathscr{T}^T, \, i = 1, \ldots, k.
\]
If we use a basis $\mathscr T^f = \spa\{ \bs t^{f}_1,\ldots, \bs t^{f}_{\ell} \}$, then for $0\leq k\leq \ell$, 
\[
\alt^k(\mathscr T^f) = {\rm span}\left \{\dd  t_{\sigma}^{f} := \dd  t^{f}_{\sigma(1)}\wedge\ldots \wedge \dd  t_{\sigma(k)}^{f} \mid \sigma \in \Sigma(k, \ell) \right \}.
\]
If $\{\bs t_1^f,\ldots,\bs t_\ell^f\}$ is orthonormal, then this basis is self-dual.
Since $\bs t^f_i \in \mathbb R^d$ and each increasing sequence \(\sigma \in \Sigma(k, \ell)\) can be naturally embedded into \(\Sigma(k, d)\), we have $\dd  t_{\sigma}^{f} \in \alt^k(\mathscr T^T)$. Therefore, we will omit the $\Pi_f^*$ operator when expressing differential forms using $\dd  t_{\sigma}^{f}$.

For a polynomial differential form \(\omega = \sum_{\sigma \in \Sigma(k, \ell)} a_\sigma \dd t_{\sigma}^{f} \in \mathbb{P}_r\Lambda^k(f)\), we can embed it into \(\mathbb{P}_r\Lambda^k(T)\) as follows:
\begin{itemize}[leftmargin=16pt]
\item Extend \(a_{\sigma}\in \mathbb{P}_r(f)\) using the Bernstein basis \(\mathbb{P}_r(f) \to \mathbb{P}_r(T)\).
\item Extend \(\dd t_{\sigma}^{f}\) using the embedding \(\alt^k(\mathscr{T}^f) \hookrightarrow \alt^k(\mathscr{T}^T)\) defined above.
\end{itemize}

\smallskip

Similarly, we can embed $\alt^k(\mathscr{N}^f) \hookrightarrow \alt^k(\mathscr{T}^T)$ and express a basis using a normal basis $\mathscr N^f = \spa\{ \bs n^{f}_1,\ldots, \bs n^{f}_{d - \ell} \}$. For $0\leq k\leq d - \ell$, we have 
\[
\alt^k(\mathscr N^f) = {\rm span}\left \{\dd  n_{\sigma}^{f}: = \dd  n^{f}_{\sigma(1)}\wedge\ldots \wedge \dd  n_{\sigma(k)}^{f}, \sigma \in \Sigma(k, d-\ell) \right \}.
\]
In our context, $\{ \bs n^{f}_1,\ldots, \bs n^{f}_{d - \ell} \}$ is usually non-orthogonal.

\subsection{Traces of differential forms}
Let $f\in\Delta_\ell(T)$, and let $\iota_f:f\hookrightarrow T$ denote the inclusion. For $0\leq k\leq d$, define the tangential trace by pullback:
$$
\tr_f:=\iota_f^*:\Lambda^k(T)\to\Lambda^k(f).
$$
Thus, for $x\in f$ and $v_1,\ldots,v_k\in\mathscr T^f$,
$$
(\tr_f\omega)(v_1,\ldots,v_k)=\omega(v_1,\ldots,v_k).
$$
If $\dim f<k$, then $\Lambda^k(f)=\{0\}$ and hence $\tr_f\omega=0$.

For $\omega=a_\sigma(x)\dd y_\sigma$, we have
$$
\tr_f\bigl(a_\sigma(x)\dd y_\sigma\bigr)=\bigl.a_\sigma(x)\bigr|_f\,\tr_f\dd y_\sigma.
$$
Therefore, $\tr_f\bigl(a_\sigma(x)\dd y_\sigma\bigr)=0$ if either $\bigl.a_\sigma(x)\bigr|_f=0$ or $\tr_f\dd y_\sigma=0$. In particular, if $\dd y_\sigma$ contains $\dd n$ for some $\bs n\in\mathscr N^f$, then $\tr_f\dd y_\sigma=0$, since $\dd n$ vanishes on $\mathscr T^f$. For example, $\tr_{F_i}\dd\lambda_f=0$ for $i\in f$, because $\dd\lambda_f$ contains $\dd\lambda_i=(\nabla\lambda_i)^\flat$ and $\nabla\lambda_i\perp F_i$.

For the $1$-form $\dd\lambda_i$,
$$
\tr_f\dd\lambda_i=\dd_f\lambda_i.
$$
The vector proxies of $\dd\lambda_i$ and $\dd_f\lambda_i$ are $\nabla\lambda_i$ and $\nabla_f\lambda_i$, respectively. Since the trace is a pullback, it preserves wedge products and commutes with the exterior derivative:
$$
\tr_f(\omega\wedge\eta)=\tr_f\omega\wedge\tr_f\eta,
\qquad
\tr_f(\dd\omega)=\dd_f(\tr_f\omega).
$$
%In particular,
%$$
%\tr_f(\dd\lambda_i\wedge\dd\lambda_j)=\dd_f\lambda_i\wedge\dd_f\lambda_j.
%$$
%where $\nabla_f = \Pi_f\nabla$ is the surface gradient on $f$.

By treating $f$ as the ambient simplex and considering a subsimplex $e \in \Delta_{\ell-1}(f)$, we can define the boundary trace:
\[
\tr_{\partial f}: \Lambda^k(f) \to \Lambda^k(\partial f), \quad (\tr_{\partial f} \omega)|_e = \tr_e \omega, \quad e \in \Delta_{\ell-1}(f), \quad 0 \leq k < \ell \leq d.
\]

A lesser-known trace is the normal trace~\cite{AwanouFabienGuzmanStern2023,KurzAuchmann2012,Weck2004,bossavit1991differential}:
% \cite[Section 4.3.1]{bossavit1991differential}:
\[
\tr_F^n: \Lambda^k(T) \mapsto \Lambda^{k-1}(F), \quad F \in \Delta_{d-1}(T), \, 1 \leq k \leq d,
\]
which is defined using the outward unit normal vector $\boldsymbol{n}_F$ of the face $F$ as 
\[
\tr_F^n \omega (v_1, \ldots, v_{k-1}) := \omega \, \lrcorner\, \, \bs n_F =\omega (\bs n_F, v_1, \ldots, v_{k-1}) |_F,
\]
for $v_i \in \mathscr{T}^F, \, i = 1, \ldots, k-1.$
If $\omega = \dd n_F\wedge \eta$ with $\eta\in \Lambda^{k-1}(F)\hookrightarrow \Lambda^{k-1}(T)$, then $\tr_F^n \omega = \eta$. 

Similarly, for $e \in \Delta_{\ell-1}(f)$, let $\bs n_f^e = - \nabla_f \lambda_{f\setminus  e}/ |\nabla_f \lambda_{f\setminus  e}|$ be the outward unit normal vector to $e$ within $\mathscr T^f$. We define
\[
\tr_e^n: \Lambda^k(f) \mapsto \Lambda^{k-1}(e), \qquad \tr_e^n \omega := \omega \, \lrcorner\, \, \bs n_f^e, \qquad 1 \leq k \leq \ell.
\]
Notice that the normal trace is defined for $k$-forms with $k \geq 1$, excluding the case $k = 0$.

When $\omega \in \Lambda^k( f)\hookrightarrow \Lambda^k(T)$, for $1\leq k\leq d-1$, we have 
\begin{equation}\label{eq:traceofextension}
 \tr_f \omega = \omega, \quad \tr_F^n\omega = 0, \qquad  \text{ for } f\in \Delta(F), F\in \Delta_{d-1}(T).
\end{equation}

We have the following duality between the tangential and normal traces; see~\cite[Section 2.2]{AwanouFabienGuzmanStern2023}. 
%and~\cite[Section 4.3.1]{bossavit1991differential} (with a correction on the sign). 

\begin{lemma}
For $\omega\in \Lambda^k(T), 0\leq k\leq d-1,$ and $F\in \Delta_{d-1}(T)$, we have 
\begin{equation}\label{eq:duality}
\star_F  \tr_F \omega = (-1)^{k} \tr_F^n \star_T \, \omega,
\end{equation}
where $\star_T$ is the Hodge star for a positively oriented $\mathscr T^T$ and $\star_F$ is the Hodge star on the $(d-1)$-dimensional space $\mathscr T^F$ with the induced orientation. Similarly,
for $e\in \Delta_{\ell-1}(f)$ and $f\in \Delta_{\ell}(T)$ with $1\leq \ell \leq d$,
\[
\star_e \tr_e \omega = (-1)^{k}  \tr_e^n \star_f \, \omega, \quad \omega\in \Lambda^k(f), 0\leq k\leq \ell - 1.
\]  
%where 
\end{lemma}
\begin{proof}
Let $\{\bs t_1,\ldots,\bs t_{d-1}\}$ be an orthonormal basis of $\mathscr T^F$. Choose the orientation of $F$ induced by the positive orientation of $T$, so that $(\bs n_F,\bs t_1,\ldots,\bs t_{d-1})$ is a positively oriented orthonormal $t$-$n$ frame of $\mathbb R^d$. Let $(x_1,\ldots,x_d)$ be the associated Cartesian coordinates. We distinguish two cases according to $\sigma\in\Sigma(k,d)$.

\step{1} $1\in\sigma$ (normal component present).\; Then $\dd x_{\sigma}$ contains $\dd n_F$ so $\operatorname{tr}_F \dd x_{\sigma}=0$.  Because $1\notin\sigma^c$, the $(d - k)$-form $\star_T  \dd x_{\sigma} = \pm \dd x_{\sigma^c}$ is purely tangential; hence \(\operatorname{tr}^n_F(\star_T \dd x_{\sigma})=0\). We conclude that \eqref{eq:duality} holds with both sides equal to zero.

\smallskip

\step{2} $1\notin\sigma$ (purely tangential). Then $\dd x_\sigma \in \Lambda^{k}( \mathscr T^F)$ and consequently, by \eqref{eq:traceofextension}, 
\[
\star_F\tr_F \dd x_\sigma = \star_F\dd x_\sigma=\sign (\sigma, \sigma^c\setminus\{1\})\dd x_{\sigma^c\setminus\{1\}}.
% = \sign (\sigma, \sigma^c(1:d-k-1))\dd x_{\sigma^c(1:d-k-1)}.
\]
We compute the right-hand side of \eqref{eq:duality}: $\star_T \dd x_\sigma = \sign (\sigma, \sigma^c) \dd x_{\sigma^c}$. As $1\in\sigma^c$, we can write $\dd x_{\sigma^c} =  \dd n_F\wedge \dd x_{\sigma^c \setminus\{1\}}$ and calculate 
$$
\tr_F^n \dd x_{\sigma^c} = \tr_F^n (\dd n_F \wedge \dd x_{\sigma^c \setminus\{1\}}) = \dd x_{\sigma^c \setminus\{1\}}.
$$ 
The sign is $\sign (\sigma, \sigma^c\setminus\{1\})\sign (\sigma, \sigma^c) = (-1)^{k}$. The proof is thus completed.
\end{proof}
%\RV{There is a consistency issue of two Hodge stars which are defined using some }
%\subsection{Integration by parts}

%\subsection{Conformity}

%We define the {\rm coderivative} $\dd^*: \Lambda^{k} \to \Lambda^{k-1}$ as
%$$
%\dd^* := (-1)^{d(k-1) + 1} \star \, \dd \, \star.
%$$
%We have the following Green's formula~\cite[Section 4.3.1]{bossavit1991differential}:
%$$
%\langle \dd \omega, \eta \rangle_{\Omega}
%=
%\langle \omega, \dd^*\eta \rangle_{\Omega}
%+
%\langle \tr_{\partial \Omega} \omega, \tr_{\partial \Omega}^n \eta \rangle_{\partial \Omega}.
%$$
%Therefore, for a piecewise smooth differential form on a triangulation $\mathcal{T}_h$, it belongs to $H\Lambda^k(\Omega)$ if and only if \(\tr_F \omega\) is continuous for every \(F \in \Delta_{d-1}(\mathcal{T}_h)\).

%\begin{lemma}[cf. Lemma 5.1 in \cite{ArnoldFalkWinthe2006Finite}]\label{lem:Arnold200651}
%Let $\omega\in L^2\Lambda^k(\Omega)$ be piecewise smooth with respect to the
%triangulation $\mathcal{T}_h$. Then $\omega\in H\Lambda^k(\Omega)$ if and only if $\tr_f\omega$ is single-valued for all $f\in\Delta_{\ell}(\mathcal{T}_h)$ with $k\leq\ell\leq d-1$.
%\end{lemma}

%Again we can apply the result to a subsimplex $f$ and obtain similar formulae %The above definition can be coupled with the evaluation of the coefficient function on the face $f$. 
%\section{Bubble Polynomial Differential Forms}
%
%\subsection{Full polynomial}
%
%
%\subsection{Trimmed polynomial}

\section{Tangential-Normal Bases}\label{sec:tnbases}
In this section, we present bases of $\alt^k(\mathscr T^T)$ using different $t$-$n$ decompositions on subsimplices. 

\subsection{Geometric decomposition of differential forms}
We will present a geometric decomposition of $\alt^k(\mathscr T^T)$ by fixing an $e\in \Delta_s(T)$ for $0\leq s\leq d$. 
Let $$\{\boldsymbol{t}_1^e, \ldots, \boldsymbol{t}_{s}^e, \boldsymbol{n}_1^e, \ldots, \boldsymbol{n}_{d-s}^e\}$$ be a $t$-$n$ basis of $\mathbb{R}^d$ associated with $e$. The tangential basis $\{\boldsymbol{t}_1^e, \ldots, \boldsymbol{t}_{s}^e\}$ can be chosen to be orthonormal. 

\begin{lemma}
For a fixed $e\in\Delta_s(T)$ with $0\leq s\leq d$, the set
\begin{equation}\label{eq:tnalt}
\begin{aligned}
\big\{\dd t_\sigma^e\wedge\dd n_\tau^e\mid\, &\sigma\in\Sigma(k-(\ell-s),s),\ \tau\in\Sigma(\ell-s,d-s),\\
&\ell=\max\{s,k\},\ldots,\min\{k+s,d\}\big \}
\end{aligned}
\end{equation}
is a basis of $\alt^k(\mathscr T^T)$.
\end{lemma}
\begin{proof}
Choose coordinates $\{y_i\}_{i=1}^d$ adapted to $e$ such that
$$
\dd y_i=\dd t_i^e\quad\text{for }1\leq i\leq s,
\qquad
\dd y_{s+j}=\dd n_j^e\quad\text{for }1\leq j\leq d-s.
$$
Then $\{\dd y_\sigma\}_{\sigma\in\Sigma(k,d)}$ is a basis of $\alt^k(\mathscr T^T)$. For $\sigma\in\Sigma(k,d)$, set
$$
\sigma_t:=\sigma\cap\{1,\ldots,s\},
\qquad
\sigma_n:=\sigma\cap\{s+1,\ldots,d\}.
$$
Let $\ell:=s+|\sigma_n|$ and $\tau:=\sigma_n-s$. Then
$$
\sigma_t\in\Sigma(k-(\ell-s),s),
\qquad
\tau\in\Sigma(\ell-s,d-s),
$$
and
$$
\dd y_\sigma=\dd y_{\sigma_t}\wedge\dd y_{\sigma_n}
=\dd t_{\sigma_t}^e\wedge\dd n_\tau^e.
$$
Since $|\sigma_n|$ ranges from $\max\{0,k-s\}$ to $\min\{k,d-s\}$, $\ell$ ranges from $\max\{s,k\}$ to $\min\{k+s,d\}$. The map $\sigma\mapsto(\sigma_t,\tau)$ is bijective; its inverse sends $(\sigma_t,\tau)$ to $\sigma_t\cup(s+\tau)$. Thus, the family in \eqref{eq:tnalt} is precisely the coordinate basis $\{\dd y_\sigma\}_{\sigma\in\Sigma(k,d)}$.
\end{proof}
Let $f\in\Delta_\ell(T)$ and $e\subseteq f$. Recall that $\star_f$ is the Hodge star on $f$, which defines an isomorphism $\star_f:\alt^{\ell-k}(\mathscr T^f)\to\alt^k(\mathscr T^f)$. We identify $\alt^k(\mathscr T^f)$ with its natural embedding in $\alt^k(\mathscr T^T)$.

\begin{theorem}\label{thm:Altdec}
For a fixed $e\in\Delta_s(T)$ with $0\leq s\leq d$ and $0\leq k\leq d$,
\begin{equation}\label{eq:decTe}
\alt^k(\mathscr T^T)=\Oplus_{\ell=\max\{s,k\}}^{\min\{k+s,d\}}\Oplus_{f\in\Delta_\ell(T),\,f\supseteq e}\star_f\alt^{\ell-k}(\mathscr T^e).
\end{equation}
\end{theorem}

\begin{proof}
For each admissible $\ell$ and $\tau\in\Sigma(\ell-s,d-s)$, define
$$
f=e\sqcup e^*(\tau)\in\Delta_\ell(T).
$$
Then $\dim f=s+|\tau|=\ell$. Apply the basis decomposition~\eqref{eq:tnalt}, choosing the normal basis $\{\boldsymbol n_1^e,\ldots,\boldsymbol n_{d-s}^e\}$ as the tangential-normal basis $\{\nabla_{e\cup\{i\}}\lambda_i\mid i\in e^*\}$. With this choice,
$$
\dd n_\tau^e=\dd_f\widehat{\lambda}_{[f\setminus e]}\in\alt^{\ell-s}(\mathscr N_f^e)\subset\alt^{\ell-s}(\mathscr T^f).
$$
Therefore,
for $\sigma\in\Sigma(k-(\ell-s),s)$,
$$
\dd t_\sigma^e\wedge\dd n_\tau^e\in\alt^k(\mathscr T^f).
$$
Let $\sigma^c\in\Sigma(\ell-k,s)$ be the complementary sequence, so that $\sigma\sqcup\sigma^c=\{1,\ldots,s\}$. By the definition of the Hodge star and the orthogonal decomposition $\mathscr T^f=\mathscr T^e\oplus\mathscr N_f^e$, there exists $c_{\sigma,\tau}^f\in\mathbb R\setminus\{0\}$ such that
$$
\star_f\dd t_{\sigma^c}^e=c_{\sigma,\tau}^f\bigl(\dd t_\sigma^e\wedge\dd n_\tau^e\bigr).
$$
Since $\sigma\mapsto\sigma^c$ is bijective, $\{\dd t_{\sigma^c}^e\mid\sigma^c\in\Sigma(\ell-k,s)\}$ is a basis of $\alt^{\ell-k}(\mathscr T^e)$. Hence,
$$
\left\{\dd t_\sigma^e\wedge\dd n_\tau^e\mid\sigma\in\Sigma(k-(\ell-s),s)\right\}
$$
is a basis of $\star_f\alt^{\ell-k}(\mathscr T^e)$.

As $\tau$ ranges over $\Sigma(\ell-s,d-s)$, $f=e\sqcup e^*(\tau)$ ranges over all $f\in\Delta_\ell(T)$ containing $e$. Regrouping the basis in~\eqref{eq:tnalt} by $(\ell,f)$ gives the direct-sum decomposition~\eqref{eq:decTe}.
\end{proof}

\begin{figure}[htbp]
\begin{center}
\includegraphics[width=3.7in]{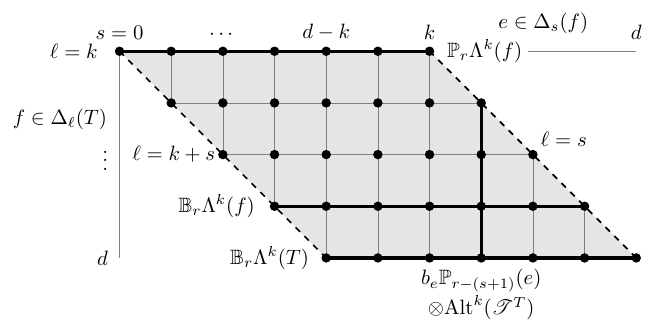}
\caption{A grid on the $(e,f)$ plane, $f\in \Delta_{\ell}(T), e\in \Delta_s(f)$, $\ell = k, \ldots, d, s = \ell -k, \ldots, \ell$, for the $t$-$n$ decomposition of polynomial differential forms $\mathbb{P}_r\Lambda^{k}(T)$ on a simplex $T$. Each dot represents the space $b_e\mathbb{P}_{r - (s+1)}(e) \otimes  \star_f \alt^{\ell - k}(\mathscr{T}^e)$. 
Summing vertically gives the geometric decomposition of $b_e\mathbb{P}_{r - (s+1)}(e) \otimes \alt^k(\mathscr{T}^T)$, and summing horizontally gives the bubble polynomial form space $\mathbb{B}_r\Lambda^k(f)$. Summing over all entries yields the full space $\mathbb{P}_r\Lambda^k(T)$.}
\label{fig:tngrid}
\end{center}
\end{figure}

In Fig.~\ref{fig:tngrid}, vertical lines inside the trapezoid illustrate such a decomposition. Due to the trapezoidal shape, the range of the summation index $\ell$ depends on $s, k,$ and $d$.

\subsection{Various $t$-$n$ bases}
% For a given $e\in \Delta_s(T)$, let $\{\boldsymbol{t}_1^e, \ldots, \boldsymbol{t}_{s}^e, \boldsymbol{n}_1^e, \ldots, \boldsymbol{n}_{d-s}^e\}$ be a  $t$-$n$ basis of $\mathbb{R}^d$ associated to $e$. And $\{\boldsymbol{t}_1^e, \ldots, \boldsymbol{t}_{s}^e\}$ can be chosen as orthonormal. 
% Then the fact
% \begin{equation*}%\label{eq:tnalt}
% \begin{aligned}
% \alt^k(\mathscr T^T) = {\rm span} \{ \dd t^e_{\sigma}\wedge \dd n^e_{\tau} \mid \, &\sigma \in \Sigma( k - (\ell - s), s), \tau \in \Sigma (\ell -s, d - s), \\
% &\ell = \max\{s,k\}, \ldots, \min\{k+s,d\} \}
% \end{aligned}
% \end{equation*}
% can be easily proved by the Vandermonde's identity \eqref{Vandermonde} as in the proof of Theorem \ref{thm:Altdec}. 
We now provide explicit bases for \eqref{eq:tnalt} via different $t$--$n$ decompositions, using a concrete choice of normal directions.
One is the face-normal basis $\{\nabla \lambda_i \mid i\in e^*\}$. By mapping $\tau \in \Sigma (\ell -s, d - s)$ to $f = e\sqcup e^*(\tau)\in \Delta_{\ell}(T)$, we can write 
\[
\dd n^e_{\tau} \to \dd \lambda_{[f\setminus  e]}.
\]
%$$
%\begin{aligned}
%\alt^k(\mathscr T^T) = {\rm span} \{ \dd t^e_{\sigma}\wedge \dd \lambda_{[f\setminus  e]} \mid \, &\sigma \in \Sigma( k - (\ell - s), s), f\in \Delta_{\ell}(T), e\subseteq f, \\
%&\ell = \max\{s,k\}, \ldots, \min\{k+s,d\} \}.
%\end{aligned}
%$$
Although $\dd\lambda_{[f\setminus e]}\in\alt^{\ell-s}(\mathscr N^e)$, it does not generally belong to $\alt^{\ell-s}(\mathscr T^f)$. We therefore replace it by $\dd_f\lambda_{[f\setminus e]}$, defined in~\eqref{eq:dflambdaf-e}. A scaled dual basis is the tangential-normal basis $\dd_f\widehat{\lambda}_{[f\setminus e]}$, defined in~\eqref{eq:dualdflambdaf-e}.

%Notice that $\dd \lambda_{[f\setminus  e]}\in \alt^{\ell -s}(\mathscr N^e)$ but not in $\alt^{\ell -s}(\mathscr T^f)$. We can change to $\dd_f \lambda_{[f\setminus  e]}$ defined in \eqref{eq:dflambdaf-e}. To be biorthogonal to $\dd \lambda_{[f\setminus  e]}$, we will use the tangential-normal basis $\dd_f \widehat{\lambda}_{[f\setminus  e]}$ defined in \eqref{eq:dualdflambdaf-e}. 

\begin{corollary}\label{cor:dualtnbasis}
Let $\{\boldsymbol{t}_1^e, \ldots, \boldsymbol{t}_{s}^e\}$ be an orthonormal basis of $\mathscr T^e$. Then
\begin{equation}\label{eq:tnbasis1}
\begin{aligned}
 \{ \dd t^e_{\sigma}\wedge \dd_f \widehat{\lambda}_{[f\setminus  e]} \mid \, &\sigma \in \Sigma( k - (\ell - s), s), f\in \Delta_{\ell}(T), f\supseteq e, \\
&\ell = \max\{s,k\}, \ldots, \min\{k+s,d\} \},
\end{aligned}
\end{equation}
and
\begin{equation*}%\label{eq:tnbasis2}
\begin{aligned}
\{ \dd t^e_{\sigma}\wedge \dd \lambda_{[f\setminus  e]} \mid \, &\sigma \in \Sigma( k - (\ell - s), s), f\in \Delta_{\ell}(T), f\supseteq e, \\
&\ell = \max\{s,k\}, \ldots, \min\{k+s,d\} \}
\end{aligned}
\end{equation*}
are scaled dual bases of $\alt^k(\mathscr T^T)$.
\end{corollary}
\begin{proof}
As $\{\boldsymbol{t}_1^e, \ldots, \boldsymbol{t}_{s}^e\}$ is orthonormal, components corresponding to different values of $\ell$ are orthogonal, since they have different normal degrees with respect to the orthogonal splitting
$\mathscr T^T=\mathscr T^e\oplus^\perp\mathscr N^e$.
Hence, for a fixed $\ell$, it suffices to prove
\begin{equation*}
\langle \dd \lambda_{[f\setminus  e]}, \dd_g \widehat{\lambda}_{[g\setminus  e]}\rangle=\delta_{f,g}c_{f\setminus  e},\quad  f,g\in\Delta_{\ell}(T),\quad e\subseteq f\cap g,\quad c_{f\setminus  e}\neq 0.
\end{equation*}
When $f\neq g$, there exists an $m\in e^*\cap f$ such that $m\not\in e^*\cap g$.
By Lemma~\ref{lem:normaldual}, $\nabla\lambda_m\perp \nabla_{e\cup \{(e^*\cap g)(i)\}} \lambda_{(e^*\cap g)(i)}\in \mathscr T^g$ for $i=1,\ldots, \ell-s$.
Then $\langle \dd \lambda_{[f\setminus  e]}, \dd_g \widehat{\lambda}_{[g\setminus  e]}\rangle=0$.

When $f=g$,
\begin{align*}
\langle \dd \lambda_{[f\setminus  e]}, \dd_f \widehat{\lambda}_{[f\setminus  e]}\rangle 
%&= \langle\dd n_{F_{(e^*\cap f)(1)}}\wedge\ldots\wedge\dd n_{F_{(e^*\cap f)(\ell-s)}}, \dd n_{f,1}^e\wedge\ldots\wedge\dd n_{f,\ell-s}^e\rangle \\
&= \det(\langle \dd \lambda_{[f\setminus  e](i)}, \dd_{e\cup \{[f\setminus  e](j)\}} \lambda_{[f\setminus  e](j)} \rangle) \\
&=\Pi_{i=1}^{\ell-s} |\nabla_{e\cup \{[f\setminus  e](i)\}} \lambda_{[f\setminus  e](i)}|^2\neq0.
\end{align*}
This completes the proof.
\end{proof}

For larger $k$, it is more convenient to express the basis \eqref{eq:tnbasis1} using the Hodge star $\star_T$.

\begin{lemma}
Let $\{\boldsymbol{t}_1^e, \ldots, \boldsymbol{t}_{s}^e\}$ be an orthonormal basis of $\mathscr T^e$. Then
\begin{equation}\label{eq:hodgestar4basis}
\star_T(\dd t_{\sigma}^e\wedge\dd_f \widehat{\lambda}_{[f\setminus  e]}) = c_{\sigma}^f(\dd t_{\sigma^c}^e\wedge\dd \lambda_{[f^*]}), \quad \sigma \in \Sigma( k - (\ell - s), s),
\end{equation}
where $c_{\sigma}^f = \frac{\left|\dd t_{\sigma}^e\wedge\dd_f \widehat{\lambda}_{[f\setminus  e]}\right|^2}{\star_T((\dd t^e_{\sigma}\wedge \dd_f \widehat{\lambda}_{[f\setminus  e]}) \wedge (\dd t_{\sigma^c}^e\wedge\dd \lambda_{[f^*]}))}$ and $\sigma^c\in \Sigma( \ell - k, s), \sigma\cup \sigma^c = \{1,\ldots,s\}$.
\end{lemma}

\begin{proof}
Set $\omega:=\dd t_\sigma^e\wedge\dd_f\widehat{\lambda}_{[f\setminus e]}$. By construction, $\dd t_\sigma^e$ is tangential to $e$, whereas $\dd_f\widehat{\lambda}_{[f\setminus e]}$ lies in the normal directions within $f\supseteq e$. Moreover, $\dd\lambda_{[f^*]}$ spans the normal directions complementary to $f$. Hence, by the definition of the Hodge star in~\eqref{eq:Hodge}, $\omega$ pairs nontrivially only with the $(d-k)$-form $\dd t_{\sigma^c}^e\wedge\dd\lambda_{[f^*]}$, whose wedge product with $\omega$ is a nonzero volume form. Therefore, there exists a scalar $c_\sigma^f$ such that~\eqref{eq:hodgestar4basis} holds.

To determine $c_\sigma^f$, wedge both sides of~\eqref{eq:hodgestar4basis} with $\omega$ and use the inner-product identity for $k$-forms in~\eqref{eq:innerproductHodge}. Since $\star_T\omega=c_\sigma^f(\dd t_{\sigma^c}^e\wedge\dd\lambda_{[f^*]})$, we obtain
$$
|\omega|^2\dd x=\omega\wedge\star_T\omega=c_\sigma^f\,\omega\wedge\dd t_{\sigma^c}^e\wedge\dd\lambda_{[f^*]}.
$$
%Since $\star_T\omega=c_\sigma^f(\dd t_{\sigma^c}^e\wedge\dd\lambda_{[f^*]})$, we obtain
%$$
%|\omega|^2\dd x=c_\sigma^f\,\omega\wedge\dd t_{\sigma^c}^e\wedge\dd\lambda_{[f^*]}.
%$$
Applying $\star_T$ to the $d$-form on the right, which extracts its scalar coefficient, gives
$$
c_\sigma^f=\frac{|\omega|^2}{\star_T\!\bigl(\omega\wedge\dd t_{\sigma^c}^e\wedge\dd\lambda_{[f^*]}\bigr)}.
$$
This is the stated formula and completes the proof.
\end{proof}

% \begin{proof}
% By the definition of the Hodge star operator \RV{in \eqref{eq:Hodge}} and the fact that $\dd t_{\sigma}^e$, $\dd_f \widehat{\lambda}_{[f\setminus  e]}$ and $\dd \lambda_{[f^*]}$ are normal each other, there exists a constant $c_{\sigma}^f$ such that
% \begin{equation*}
% \star_T(\dd t_{\sigma}^e\wedge\dd_f \widehat{\lambda}_{[f\setminus  e]}) = c_{\sigma}^f(\dd t_{\sigma^c}^e\wedge\dd \lambda_{[f^*]}).
% \end{equation*}
% \RV{To determine $c_{\sigma}^f$}, applying the definition of the inner product for $k$-forms \RV{in \eqref{eq:innerproductHodge}}, we compute the wedge product as
% \begin{align*}
% \left|\dd t_{\sigma}^e\wedge\dd_f \widehat{\lambda}_{[f\setminus  e]}\right|^2\dx &= \dd t_{\sigma}^e\wedge\dd_f \widehat{\lambda}_{[f\setminus  e]}\wedge\star_T(\dd t_{\sigma}^e\wedge\dd_f \widehat{\lambda}_{[f\setminus  e]})  \\
% &= c_{\sigma}^f\dd t_{\sigma}^e\wedge\dd_f \widehat{\lambda}_{[f\setminus  e]}\wedge \dd t_{\sigma^c}^e\wedge\dd \lambda_{[f^*]}.
% \end{align*}
% This completes the proof.
% \end{proof}

\begin{corollary}
For $0\leq s\leq d$ and $e\in\Delta_s(T)$, choose an orthonormal tangential frame on $e$. Then
we have % the set
\begin{equation}\label{eq:tnbasis}
\begin{aligned}
\alt^k(\mathscr T^T)  = {\rm span} \{ \star_T(\dd t_{\sigma^c}^e\wedge\dd \lambda_{[f^*]} ) \mid  \, & \sigma^c\in\Sigma(\ell-k, s), f\in\Delta_{\ell}(T), f\supseteq e, \\
&\ell=\max\{s,k\}, \ldots, \min\{k+s,d\} \},
\end{aligned}
\end{equation}
and a scaled dual basis is
\begin{align*}
% \alt^k(\mathscr T^T)  = {\rm span} 
\{ \star_T(\dd t_{\sigma^c}^e\wedge\mathop{\textstyle\bigwedge}\nolimits_{i\in[f^*]}\dd_{e\cup\{i\}}\lambda_i ) \mid  \, & \sigma^c\in\Sigma(\ell-k, s), f\in\Delta_{\ell}(T), f\supseteq e, \\
&\ell=\max\{s,k\}, \ldots, \min\{k+s,d\} \}.
\end{align*}
The wedge factors are ordered increasingly and are based on $e$. Their scaled duality follows from $\langle\dd\lambda_i,\dd_{e\cup\{j\}}\lambda_j\rangle=\delta_{ij}|\nabla_{e\cup\{i\}}\lambda_i|^2$, the orthonormality of the tangential frame, and the isometry of $\star_T$.
\end{corollary}

%\RV{modify to the new notation}.

\subsection{Examples}
When $s=0$, $e$ is a vertex. Without loss of generality, take $e = \{ 0 \}$. Since $s=0$, there are no tangential vectors on $e$, and only faces $f\in\Delta_k^0(T)$ occur. The tangential-normal vector can be chosen as the edge vector $\bs t_{0,i}$ for $i\in f(1:k)$. The basis is
\[
\alt^k(\mathscr T^{T})
=
\spa \left\{
\dd t_{0,f(1)} \wedge \dd t_{0,f(2)} \wedge \cdots \wedge \dd t_{0,f(k)}
\;\middle|\;
f\in\Delta_k^0(T)
\right\}.
\]

When $s=d$, $e = T$. Then only $f= e = T$ appears, and there are no normal components. We have
\[ 
\alt^k(\mathscr T^{T}) = \spa \{ \dd x_{\sigma}, \sigma\in\Sigma(k, d)\},
\]
which is the standard basis of $\alt^k$ using the ambient orthonormal basis $\{ \dd x_1, \ldots, \dd x_d\}$. 

When $k = 1$, for $s=0,\ldots, d$, we have
\begin{equation*}
\alt^{1}(\mathscr T^T)  = {\rm span} \{\dd t_{1}^e,\ldots, \dd t_{s}^e, \dd_{e\cup e^*(0)} \lambda_{e^*(0)},\ldots,\dd_{e\cup e^*(d-s-1)} \lambda_{e^*(d-s-1)}\}.
\end{equation*}

When $k = d-1$, for $s=0,\ldots, d$, we use the basis \eqref{eq:tnbasis} with the Hodge star:
\begin{equation*}
\alt^{d-1}(\mathscr T^T)  = {\rm span} \{ \star_T(\dd \lambda_{e^*(0)}),\ldots,\star_T(\dd \lambda_{e^*(d-s-1)}), \star_T(\dd t_{1}^e),\ldots, \star_T(\dd t_{s}^e)\}.
\end{equation*}

\begin{example}\rm
We present examples for $d = 3$ with $k = 1$ and $k = 2$ in Table~\ref{fig:tnbasesd3k1} and Table~\ref{fig:tnbasesd3k2}, respectively, corresponding to edge and face elements in three dimensions.

\begin{table}[htp]
\caption{Example: $t$-$n$ bases \eqref{eq:tnbasis1} of $(\alt^k(\mathscr T^{T}))^{\sharp}$ for $d=3$ and $k=1$: edge element.}
\begin{center}
\begin{tabular}{@{}c c c c@{}}
\toprule
$s=0$ & $s = 1$ & $s = 2$ & $s = 3$\\
\midrule
\includegraphics*[width=2.6cm]{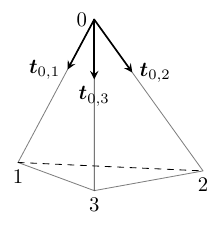} & \includegraphics*[width=2.65cm]{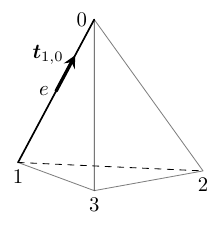}  & \includegraphics*[width=2.65cm]{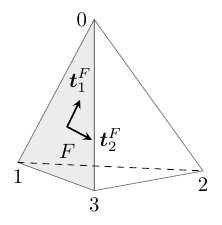}  & \\
$\ell = 1$ & $\ell = 1$ & $\ell = 2$& \\
 & \includegraphics*[width=2.75cm]{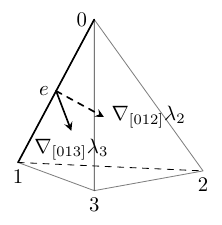} & \includegraphics*[width=2.65cm]{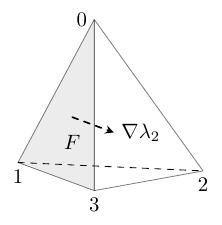} &\includegraphics*[width=2.65cm]{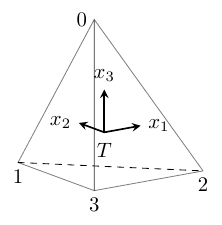}\\
 & $\ell = 2$ & $\ell = 3$& $\ell = 3$\\
\bottomrule
\end{tabular}
\end{center}
\label{fig:tnbasesd3k1}
\end{table}%

\begin{table}[htp]
\caption{Example: $t$-$n$ bases \eqref{eq:tnbasis} of $-(\star\alt^k(\mathscr T^{T}))^{\sharp}$ for $d=3$ and $k=2$: face element.}
\begin{center}
\begin{tabular}{@{}c c c c@{}}
\toprule
$s=0$ & $s = 1$ & $s = 2$ & $s = 3$\\
\midrule
\includegraphics*[width=2.65cm]{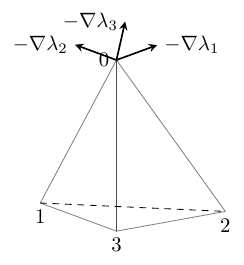} & \includegraphics*[width=3.05cm]{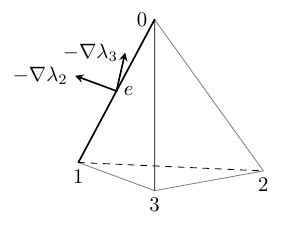}  & \includegraphics*[width=2.65cm]{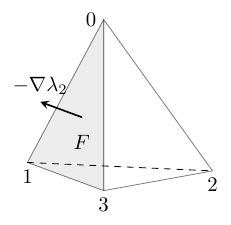}  & \\
$\ell = 2$ & $\ell = 2$ & $\ell = 2$& \\
 & \includegraphics*[width=2.75cm]{figures/3DNDs1l1} & \includegraphics*[width=2.65cm]{figures/3DNDs2l2} & \includegraphics*[width=2.65cm]{figures/3DNDs3}\\
 & $\ell = 3$ & $\ell = 3$& $\ell = 3$ \\
\bottomrule
\end{tabular}
\end{center}
\label{fig:tnbasesd3k2}
\end{table}%

\end{example}

\section{Geometric Decomposition of Full Polynomial Differential Forms}\label{sec:geodecfullpoly}
In this section, we present a geometric decomposition of the second family of polynomial differential forms $\mathbb P_r\Lambda^{k}(T) := \mathbb P_r(T)\otimes \alt^k(\mathscr T^T)$. 
%on a simplex for which the shape function space is the full polynomial space. 
%We begin with a tensor product with the Lagrange element and then choose appropriate  $t$-$n$ bases to represent $\alt^k$. We give a precise characterization of the bubble polynomial forms. Together with the geometric decomposition of Lagrange element, we can associate for each subsimplex $e$ a frame to express $\alt^k$ so that the well documented Lagrange basis can be used. 

\subsection{Lagrange element}
For a $0$-dimensional face, i.e., a vertex $\texttt{v}$, we understand that
$\int_{\texttt{v}} u \dd s = u(\texttt{v})$ for $1\in \mathbb P_r(\texttt{v})=\mathbb R$.

\begin{lemma}[Geometric decomposition of the Lagrange element, (2.6) in~\cite{ArnoldFalkWinther2009} and Theorem 2.5 in \cite{ChenHuang2024}]\label{lm:Lagrangedec}
For the polynomial space $\mathbb P_r(T)$ with $r\geq 1$ on a $d$-dimensional simplex $T$, we have the following decomposition of the shape function space
\begin{align}
\label{eq:Prdec}
\mathbb P_r(T) &= \Oplus_{\ell = 0}^d\Oplus_{f\in \Delta_{\ell}(T)} b_f\mathbb P_{r - (\ell +1)} (f).
\end{align}
%and the decomposition of DoFs
%\begin{align*}
%% \label{eq:Prstardec}
%\mathbb P_r^*(T) &= \Oplus_{\ell = 0}^d\Oplus_{f\in \Delta_{\ell}(T)} \mathcal N(\mathbb P_{r - (\ell +1)} (f)).
%\end{align*}
A function $u\in \mathbb P_r(T)$ is uniquely determined by the DoFs
\begin{equation*}%\label{eq:dofPrdec}
\int_f u p \dd s, \quad p\in \mathbb P_{r - (\ell +1)} (f), f\in \Delta_{\ell}(T), \ell = 0, \ldots, d. 
\end{equation*}
\end{lemma}
Recall that when $r < \ell + 1$, $\mathbb{P}_{r-(\ell +1)}(f) = \{0\}$.  Consequently, the last non-zero term in the decomposition \eqref{eq:Prdec} occurs for $\ell \leq \min\{r-1, d\}$. For example, with quadratic polynomials, only vertex bubbles and edge bubbles are included in the summation, while face bubbles in higher dimensions are excluded. However, the full summation notation $\Oplus_{\ell=0}^d$ is maintained for simplicity, with the implicit understanding that the range of non-zero subspaces will naturally limit the summation.

Employing the tensor product, we obtain a straightforward decomposition for the tensor-valued Lagrange element:
\begin{align}
 \notag 
\mathbb P_r\Lambda^{k}(T) &:= \mathbb P_r(T)\otimes \alt^k(\mathscr T^T) \\
 \label{eq:polynomialdecomp} 
&=  \Oplus_{s=0}^d\Oplus_{e\in\Delta_{s}(T)}\left (b_e\mathbb P_{r-(s +1)}(e)\otimes \alt^k(\mathscr T^T)\right ).
\end{align}
The natural extension of the DoFs will give the tensor-valued Lagrange element with more than the necessary continuity to be in $H\Lambda^k(\Omega)$. 

\subsection{Geometric decomposition of polynomial differential forms}
We shall further expand $\alt^k(\mathscr T^T)$ in \eqref{eq:polynomialdecomp} using decomposition \eqref{eq:decTe}, which leads to a geometric decomposition with the required continuity.
%Combining the geometric decomposition of Lagrange element and the dual tangential-normal bases \eqref{eq:DoFstnbasis1}-\eqref{eq:dualtnbasis} will induce the following geometric decomposition of $\mathbb P_r\Lambda^{k}(T)$ and $(\mathbb P_r\Lambda^{k}(T))^*$.
\begin{theorem}\label{thm:geodec}
Let $0\leq k\leq d$ and $r\geq1$. On a $d$-dimensional simplex $T$, we have the geometric decomposition
\begin{align}\label{eq:geodec1}
\mathbb P_r\Lambda^{k}(T) &= \Oplus_{s=0}^d\Oplus_{e\in\Delta_{s}(T)}\Oplus_{\ell=\max\{s,k\}}^{\min\{k+s,d\}}\\
\notag
&\quad \;\Oplus_{f\in\Delta_{\ell}(T), f\supseteq e}\Oplus_{\sigma\in\Sigma(s+k-\ell,s)}\left(b_e\mathbb P_{r-(s +1)}(e)\dd t_{\sigma}^e\wedge \dd \lambda_{[f\setminus  e]}\right).
%\\
%\label{eq:geodec2}
%(\mathbb P_r\Lambda^{k}(T))^* &= \Oplus_{s=0}^d\Oplus_{e\in\Delta_{s}(T)}\Oplus_{\ell=\max\{s,k\}}^{\min\{k+s,d\}}\\
%\notag
%&\quad \;\Oplus_{f\in\Delta_{\ell}(T), e\subseteq f}\Oplus_{\sigma\in\Sigma(s+k-\ell,s)}\mathcal N\left(\mathbb P_{r-(s +1)}(e)\dd t_{\sigma}^e \wedge \dd \lambda_{[f\setminus  e]}\right).
\end{align}
An element $\omega\in\mathbb P_r\Lambda^{k}(T)$ is uniquely determined by the DoFs:
\begin{equation}\label{eq:tndofs}
\begin{aligned}
\int_e \bigl\langle \omega,\,
\dd t_{\sigma}^e \wedge \dd_f \widehat{\lambda}_{[f\setminus e]} \bigr\rangle\,
q \,\dd s,
&\ \
e \in \Delta_s(T),\; 0\leq s\leq d,\\ 
&\ \ f \in \Delta_{\ell}(T), f\supseteq e, \max\{s,k\}\leq \ell \leq \min\{k+s,d\}, \\
&\;\;
q \in \mathbb P_{r-(s+1)}(e),\;
\sigma \in \Sigma(s+k-\ell, s).
\end{aligned}
\end{equation}
\end{theorem}
\begin{proof}
Substituting the $t$--$n$ decomposition of Theorem~\ref{thm:Altdec} into \eqref{eq:polynomialdecomp}, using the explicit basis in Corollary~\ref{cor:dualtnbasis}, gives the direct-sum decomposition~\eqref{eq:geodec1}.
We prove that each element of $\mathbb P_r\Lambda^{k}(T)$ is uniquely determined by the DoFs \eqref{eq:tndofs}.
By the geometric decomposition \eqref{eq:geodec1}, the number of degrees of
freedom in \eqref{eq:tndofs} matches the dimension of
$\mathbb P_r \Lambda^{k}(T)$.

% Assume $\omega\in\mathbb P_r\Lambda^{k}(T)$ satisfies all the DoFs \eqref{eq:tndofs} vanish.
Let $e'\in\Delta_{s'}(T)$ denote the scalar-bubble carrier of a basis function
$$
\phi_{e'}=b_{e'}p\,\dd t_{\sigma'}^{e'}\wedge \dd \lambda_{[f'\setminus e']},
\qquad p\in\mathbb P_{r-(s'+1)}(e').
$$
Since $b_{e'}|_e=0$ whenever $e'\not\subseteq e$, the DoF--basis matrix \eqref{eq:lowertriangular} is block lower triangular with respect to the dimensions of the subsimplices:
\begin{equation}\label{eq:lowertriangular}
\renewcommand{\arraystretch}{1.15}
\begin{array}{c|*{5}{c}}
\displaystyle \int_e \,\cdot \;\backslash\; \phi_{e'}
& 0 & 1 & \cdots & d-1 & d \\ \hline
0
& \graysquare & 0 & \cdots & 0 & 0 \\
1
& \graysquare & \graysquare & \cdots & 0 & 0 \\
\vdots
& \vdots & \vdots & \ddots & \vdots & \vdots \\
d-1
& \graysquare & \graysquare & \cdots & \graysquare & 0 \\
d
& \graysquare & \graysquare & \cdots & \graysquare & \graysquare
\end{array}
\end{equation}

%When restricted to the diagonal block of the DoF-Basis matrix on each subsimplex, it is block diagonal and invertible by the duality of the basis, cf. Corollary~\ref{cor:dualtnbasis} and Lemma~\ref{lm:Lagrangedec}.
For a diagonal block with $e'=e$, Corollary~\ref{cor:dualtnbasis} makes the exterior-frame pairing diagonal in $(\ell,f,\sigma)$. For each diagonal frame entry, the scalar block is the weighted mass matrix
$$
\left(\int_e b_e p_i p_j\,\dd s\right)_{i,j},\qquad\{p_i\}\text{ a basis of }\mathbb P_{r-s-1}(e).
$$
This matrix is positive definite, and hence invertible, because $b_e>0$ in the interior of $e$. The block-triangular matrix is therefore invertible, and the DoFs are unisolvent.
%
%\RV{Use the lexicographic  ordering in $(s, \ell)$ grid. The trace operator $\tr_e \omega$ will satisfy the property $\tr_{e(s', \ell')} \omega(s, \ell) = 0$ if $(s',\ell')<(s, \ell)$. So the DoF-Basis matrix is lower triangular. When restrict to the diagonal of the DoF-Basis matrix, it is diagonal by the duality of the basis.}
\end{proof}

%\RV{The geometric decomposition \eqref{eq:geodec1} differs from \eqref{eq:bubbledec}. Thanks to Theorem~\ref{th:tnbasis}, decomposition \eqref{eq:bubbledec} is used for defining the DoFs in \eqref{eq:tndofs}, whereas decomposition \eqref{eq:geodec1} is employed for constructing the basis functions.}

%\RV{In Fig.~\ref{fig:tngrid}, vertical lines inside the trapezoid illustrate such a decomposition. The range of the summation index $\ell$ depends on $s, k,$ and $d$.} \RV{Note that in \eqref{eq:tndofs} the ordering of $f$ and $e$ is reversed.
%The corresponding change in the lower and upper bounds of the summation
%is best understood from Fig.~\ref{fig:tngrid}, read rowwise.}

% \RV{modify the result using the new notation}
For larger $k$, it is more convenient to express the basis using the Hodge star.
\begin{corollary}
An element $\omega\in\mathbb P_r\Lambda^{k}(T)$ is uniquely determined by the following DoFs, for all $e\in \Delta_{s}(T), s = 0,1,\ldots, d$, and $f\in\Delta_{\ell}(T), f\supseteq e, \max\{s,k\}\leq\ell\leq\min\{k+s,d\}$:
\begin{equation}\label{eq:tndofslargek}
\begin{aligned}
\int_e \star_T(\omega \wedge \dd t_{\sigma^c}^e\wedge\dd \lambda_{[f^*]})\ q&\dd s,\quad  q\in\mathbb P_{r-(s +1)}(e), \sigma^c\in\Sigma(\ell-k, s).
\end{aligned}
\end{equation} 
\end{corollary}
\begin{proof}
By \eqref{eq:hodgestar4basis},
\begin{equation*}
\langle \omega, \dd t_{\sigma}^e\wedge \dd_f \widehat{\lambda}_{[f\setminus  e]} \rangle\dx = \omega \wedge \star_T(\dd t_{\sigma}^e\wedge\dd_f \widehat{\lambda}_{[f\setminus  e]}) = c_{\sigma}^f(\omega \wedge \dd t_{\sigma^c}^e\wedge\dd \lambda_{[f^*]}).
\end{equation*}
Thus, the DoFs \eqref{eq:tndofs} are equivalent to those in \eqref{eq:tndofslargek}.
\end{proof}

Since the DoF--basis matrix \eqref{eq:lowertriangular} is not diagonal, the basis functions induced by the decomposition \eqref{eq:geodec1} are not dual to the DoFs defined in \eqref{eq:tndofs}. To address this, in Subsection~\ref{sec:dualbasis}, we will construct basis functions dual to nodal DoFs using the Lagrange basis.

\subsection{Geometric decomposition using bubble polynomial differential forms}\label{sec:geodecBubbleform}
In this subsection, we present a geometric decomposition of $\mathbb P_r\Lambda^{k}(T)$ using bubble polynomial differential forms, which corresponds to the row sum of the subspaces in Fig.~\ref{fig:tngrid}. 
Again, we work on $T$ and obtain the result on $f$ by changing notation. 

%\subsection{Bubble polynomial differential forms}
%We will rearrange the terms in the geometric decomposition to obtain bubble differential forms. 

\begin{lemma}\label{lm:bubblefull}
For $0\leq k\leq d$ and $f\in \Delta_{\ell}(T)$ with $d-k\leq \ell \leq d$, the subspace $b_f\star_T \alt^{d-k}(\mathscr T^f)\subset \mathbb P_{\ell +1}\Lambda^k(T)$ satisfies
\begin{equation}\label{eq:bubblefull}
\tr_{\partial T} \left (b_f\star_T \alt^{d-k}(\mathscr T^f)\right ) = 0.
\end{equation}
\end{lemma}
\begin{proof}
If $k=d$, then $\tr_{\partial T}\Lambda^d(T)=0$, and the result follows. Hence, assume $0\leq k\leq d-1$. Pick an arbitrary $\omega\in b_f\star_T\alt^{d-k}(\mathscr T^f)$ and $F\in\Delta_{d-1}(T)$.

\smallskip
\step{1} Suppose $f\nsubseteq F$, or equivalently, $f\cap F^*\neq\varnothing$. Then $b_f|_F=0$, and hence $\tr_F\omega=0$. The most familiar case is $k=0$. The condition $d-k=d\leq\ell$ forces $\ell=d$ and $f=T$; thus, $b_T|_{\partial T}=0$.
\medskip

\step{2} If $f\subseteq F$, then $\ell\leq d-1$, which together with $\ell\geq d-k$ implies $1\leq k\leq d-1$. We have $\star_T \, \omega \in b_f \alt^{d-k}(\mathscr T^f)\hookrightarrow  b_f \alt^{d-k}(\mathscr T^F)$. So $\tr_F^n \star_T\, \omega = 0$. Hence, $\tr_F\omega=0$ by the identity \eqref{eq:duality}.
% \smallskip

Alternatively, we can prove the result using an explicit basis. Choose an orthonormal tangential frame $\{\bs t_1^f,\ldots,\bs t_\ell^f\}$. Since $\dd\lambda_{[f^*]}$ spans $\alt^{d-\ell}(\mathscr N^f)$, the orthogonal decomposition $\mathscr T^T=\mathscr T^f\oplus\mathscr N^f$ gives the following basis of $b_f\star_T\alt^{d-k}(\mathscr T^f)$:
\[
\left\{b_f \dd t_\sigma^f\wedge\dd\lambda_{[f^*]}\mid\sigma\in\Sigma(\ell-(d-k),\ell)\right\}.
\]
For any facet $F=F_i$, either $i\in f$ or $i\in f^*$. In the first case, $b_f$ contains the factor $\lambda_i$, so $b_f|_F=0$. In the second case, $\dd\lambda_{[f^*]}$ contains the factor $\dd\lambda_i$, whose tangential trace on $F_i$ vanishes. Thus every basis form has zero tangential trace on every facet, proving the result.
\end{proof}

%We generalize the bubble to $e\subset f$ as $b_{e} \otimes \star_{f} \alt^{\ell-k}(\mathscr T^e)$.
%\RV{not sure if it is correct. but formally I guess it is true.}
%\begin{lemma}
%Let $f\in \Delta_{\ell}(T), \ell \geq k$, and $e\in \Delta_s(f)$. For $\tilde f\in \Delta(T)$ satisfying $f\cap \tilde f^*\neq \varnothing$, then
%$$
%\tr_{\tilde f}(b_{e}\otimes \star_{f} \alt^{\ell-k}(\mathscr T^e))=0.
%$$
%\end{lemma}
%\RV{It's not true. For example, $d=3, \ell=2, k=s=1$, $e=e_{01}$, $f=F_3$, $\tilde f=F_2$. Then
%\begin{align*}
%\tr_{\tilde f}(b_{e}\otimes \star_{f} \alt^{\ell-k}(\mathscr T^e))&=\tr_{F_2}(\lambda_0\lambda_1\otimes \star_{F_3} \alt^{1}(\mathscr T^e))\\
%&=\textrm{span}\{\tr_{F_2}(\lambda_0\lambda_1\dd n_{F_3,e})\} \\
%&=\textrm{span}\{(\lambda_0\lambda_1)|_{F_2}\tr_{F_2}(\dd n_{F_3,e})\}\neq0.
%\end{align*}
%}
%\begin{proof}
%%As $f\cap \tilde f^*\neq \varnothing$, there exists a $j\in f$ but $j\not\in\tilde f$. 
%%If $j\in e$, then $b_{e}|_{\tilde f} = 0$. If $j\in e$, then  
%If $\tilde f^*\cap e\neq \varnothing$, then $b_{e}|_{\tilde f} \, = 0$.
%\end{proof}

%\begin{figure}[htbp]
%\begin{center}
%\includegraphics[width=6cm]{figures/stacked_rectangles.pdf}
%\caption{will add more explanation on this figure}
%\label{fig:dec1}
%\end{center}
%\end{figure}

Define the bubble polynomial space, for $0\leq k\leq d$ and $r\geq 1$,
\begin{align}\label{eq:bubbledecomp}    
\mathbb{B}_r \Lambda^{k}(T) &= \Oplus_{\ell=d-k}^d \Oplus_{f \in \Delta_{\ell}(T)}\left [\mathbb P_{r-(\ell +1)}(f)\otimes b_f\star_T \alt^{d-k}(\mathscr T^f) \right ]\\
& = \Oplus_{\ell=d-k}^d \Oplus_{f \in \Delta_{\ell}(T)}b_f \star_T\mathbb P_{r-(\ell +1)}\Lambda^{d-k}(f). \notag
\end{align}
%with 
%$$
%\dim\mathbb{B}_r\Lambda^{k}(T)=\sum_{\ell=d-k}^d{d+1\choose\ell+1} {r-1\choose\ell}{\ell\choose d-k}.
%$$
%\LC{Indeed, $\dd\lambda_{[f^*]}$ spans $\alt^{d-\ell}(\mathscr N^f)$, and the Hodge star maps each tangential $(d-k)$-form $\dd t_{\sigma^c}^f$ to a nonzero multiple of $\dd t_\sigma^f\wedge\dd\lambda_{[f^*]}$, where $\sigma^c$ is the complement of $\sigma$ in $\{1,\ldots,\ell\}$. Empty wedge products are understood as one.}
Then by Lemma \ref{lm:bubblefull}, for $0\leq k\leq d$, 
\begin{equation*}%\label{eq:bubbleformk0}
\mathbb{B}_r \Lambda^{k}(T)\subseteq (\mathbb{P}_r \Lambda^{k}(T)\cap \ker(\tr_{\partial T})).
\end{equation*}
We will later prove in Corollary~\ref{eq:Lambdakbubblecharac} that this inclusion is in fact an equality.
%For $k=d$, by the geometric decomposition of Lagrange element,
%\[
%\mathbb{B}_r \Lambda^{d}(T) = \mathbb P_{r}\Lambda^d(T). 
%\]
%It also satisfies $\tr_{\partial T}\mathbb{B}_r \Lambda^{d}(T) = 0$ as no non-trivial $d$-form on $(d-1)$-dimensional faces.
For the lowest-order top-degree endpoint, we adopt the convention
$\mathbb B_0\Lambda^d(T):=\mathbb P_0\Lambda^d(T).$
%\begin{remark}\rm 
%%We shall show later the equality holds. 
% We have the following relation:
%$$
%\mathbb B_{r}\Lambda^k(T) \subset \mathbb D_{r}\Lambda^k(T) \cong \mathbb B_{r+1}^{-}\Lambda^k(T)\subseteq \mathbb B_{r+1}\Lambda^k(T).
%$$
%
%%Recall that the isomorphism in the middle is given by the mapping $\dd \lambda_{[f]} \leftrightarrow w_f$.
%\end{remark}

\begin{example}\rm
For $k=1$, \eqref{eq:bubbledecomp} becomes
\begin{equation*}
(\mathbb{B}_r \Lambda^{1}(T))^{\sharp}=\Oplus_{F \in \Delta_{d-1}(T)}\left [\mathbb P_{r-d}(F) \otimes b_F\mathscr N^F\right ] \Oplus (\mathbb P_{r-(d+1)}(T)\otimes b_T\mathbb R^d).
\end{equation*}
This characterization of the edge element bubble space was first given in \cite{ChenChenHuangWei2024}. 

Let us then consider the case $k = d-1$:
\begin{equation*}  
(\star_T\mathbb{B}_r \Lambda^{d-1}(T))^{\sharp}=\Oplus_{\ell=1}^d \Oplus_{f \in \Delta_{\ell}(T)}\left [\mathbb P_{r-(\ell +1)}(f) \otimes b_f\mathscr T^f \right ].
\end{equation*}
This characterization of the $H(\div)$ bubble finite element space was first given in \cite{ChenHuang2024}. 
% $\Box$
\end{example}

In general, for a subsimplex $f\in \Delta_{\ell}(T)$ with $\ell \geq k$, define
\begin{align}\label{eq:Brf}
\mathbb{B}_r \Lambda^{k}(f) &:= \Oplus_{s=\ell - k}^{\ell} \Oplus_{e \in \Delta_{s}(f)}\left [\mathbb{P}_{r-(s+1)}(e)\otimes b_e\star_{f} \alt^{\ell - k}(\mathscr T^e) \right ]\\
&\;= \Oplus_{s=\ell - k}^{\ell} \Oplus_{e \in \Delta_{s}(f)} b_e \star_{f} \mathbb{P}_{r-(s+1)}\Lambda^{\ell - k}(e). \notag
\end{align}
That is, for a fixed $f$, we sum the subspaces horizontally in Fig.~\ref{fig:tngrid}.

When $\ell = \dim f \geq k$, we know $\tr_{\partial f} \mathbb{B}_r \Lambda^{k}(f) = 0$ by applying Lemma \ref{lm:bubblefull} with the changes in notation $T\to f$, $f\to e$ and $d\to \ell$. 

When $\ell = \dim f = k$, by the geometric decomposition of the Lagrange element,
\begin{align*}
\mathbb{B}_r \Lambda^{k}(f) = &\star_{f}  \Oplus_{s=0}^{k} \Oplus_{e \in \Delta_{s}(f)}b_e \mathbb{P}_{r-(s+1)}\Lambda^{0}(e) = \star_{f} \mathbb P_r\Lambda^0(f) = \mathbb P_r\Lambda^k(f).
\end{align*}
It also satisfies $\tr_{\partial f}\mathbb{B}_r \Lambda^{k}(f) = 0$ because there are no nontrivial $k$-forms on $(k-1)$-dimensional subsimplices.

\begin{remark}\label{rm:bubblespace}\rm
The consistent extension space $E_f\mathring{X}(f)$ defined in~\cite{ArnoldFalkWinther2009} has vanishing tangential trace on every face of $T$ that does not contain $f$; see~\cite[Lemma~4.2]{ArnoldFalkWinther2009}. By contrast, although $\mathbb B_r\Lambda^k(f)$ defined in~\eqref{eq:Brf} is intrinsically trace-free on $\partial f$, its natural embedding into $\mathbb P_r\Lambda^k(T)$ need not have vanishing tangential trace on such faces.

For example, let $d=3$, $s=k=1$, $\ell=2$, and $r=2$. Take $f=F\in\Delta_2(T)$ and $e\in\Delta_1(F)$, and let $i$ be the unique vertex in $F\setminus e$. Then $e=F\cap F_i$. With the chosen orientation, the form $b_e\star_F\dd t^e$ belongs to $\mathbb B_r\Lambda^1(F)$ and satisfies
\begin{equation*}
\tr_{F_i}\left(b_e\star_F\dd t^e\right)
=\tr_{F_i}\left(-b_e\frac{(\nabla_F\lambda_i)^\flat}{|\nabla_F\lambda_i|}\right)\neq0
\qquad\text{if }F\not\perp F_i.
\end{equation*}
Thus, in general, $\mathbb B_r\Lambda^k(f)\neq E_f\mathring{X}(f)$. Decomposition~\eqref{eq:Brf} and Isaac's decomposition~\cite[(20)]{Isaac2022} describe the same trace-free polynomial space on $f$, although their individual summands need not coincide.
\end{remark}
%\RV{please check the example. for a 1-form, only one d. Also $b_e$ is invovled. Use t-n vectors. is it $b_e n^e_F$? vanishing trace on $F$ but not another $F'$ containng $e$.}

%\subsection{Geometric decomposition}
\begin{theorem}[Geometric decomposition of the second family of polynomial differential forms]
For $0\leq k\leq d$ and $r\geq 1$, we have the geometric decomposition
 \begin{align}\label{eq:bubbledec}
\mathbb P_r\Lambda^{k}(T) &= \Oplus_{\ell=k}^d\Oplus_{f\in\Delta_{\ell}(T)}\mathbb B_r\Lambda^k(f)
\\
\notag & = \Oplus_{f\in \Delta_k(T)}\mathbb P_r\Lambda^k (f) \oplus \Oplus_{\ell=k+1}^d\Oplus_{f\in\Delta_{\ell}(T)}\mathbb B_r\Lambda^k(f),
\end{align}
where $\mathbb B_r\Lambda^k(f)$ is defined in \eqref{eq:Brf}.
\end{theorem}
\begin{proof}
%By \eqref{eq:Prf}, we only prove \eqref{eq:keydec}. 
By the definition \eqref{eq:Brf}, we have
\begin{align*}    
&\quad \;\Oplus_{\ell=k}^d\Oplus_{f\in\Delta_{\ell}(T)}\mathbb B_r\Lambda^k(f)\\
&=\Oplus_{\ell=k}^d\Oplus_{f\in\Delta_{\ell}(T)}\Oplus_{s=\ell-k}^{\ell} \Oplus_{e \in \Delta_{s}(f)}\left [\mathbb{P}_{r-(s+1)}(e)\otimes b_e\star_{f} \alt^{\ell-k}(\mathscr T^e) \right ],
\end{align*}
which is the row sum of the $t$-$n$ decomposition in Fig.~\ref{fig:tngrid}. 

Switch the order of indices $\ell$ and $s$, i.e., sum column-wise, to get
\begin{align*}    
&\quad \;\Oplus_{\ell=k}^d\Oplus_{f\in\Delta_{\ell}(T)}\mathbb B_r\Lambda^k(f)\\
&=\Oplus_{s=0}^d\Oplus_{e \in \Delta_{s}(T)}\mathbb{P}_{r-(s+1)}(e)b_e\otimes  \left [\Oplus_{\ell=\max\{s,k\}}^{\min\{k+s,d\}}\Oplus_{f\in\Delta_{\ell}(T), e\subseteq f} \star_{f} \alt^{\ell-k}(\mathscr T^e) \right ],
\end{align*}
which, together with the decomposition \eqref{eq:decTe} of $\alt^k(\mathscr T^T)$ and the geometric decomposition \eqref{eq:polynomialdecomp} of the Lagrange element, completes the proof.
%For $\ell \geq k+1$, similar to the proof of Lemma \ref{lm:bubble}, we have
%\begin{equation}
%\tr_{\partial f} \mathbb{B}_r \Lambda^{k}(f) = 0.
%\end{equation}
\end{proof}

As explained in Remark~\ref{rm:bubblespace}, the geometric decomposition \eqref{eq:bubbledec} is different
from the geometric decomposition (4.9) in \cite{ArnoldFalkWinther2009}. The geometric decomposition \eqref{eq:bubbledec} also differs from \eqref{eq:geodec1}.

%\subsection{Properties of bubble polynomial forms}
Next, we show that $\mathbb{B}_r \Lambda^{k}(T)$ contains all bubble polynomial forms.
\begin{corollary}\label{eq:Lambdakbubblecharac}
For $0\leq k\leq d-1$ and $r\geq1$,
\begin{equation*}
 \mathbb{B}_r \Lambda^{k}(T) =  (\mathbb{P}_r \Lambda^{k}(T)\cap \ker(\tr_{\partial T})).
\end{equation*}
\end{corollary}
\begin{proof}
By Lemma \ref{lm:bubblefull}, for $0\leq k\leq d-1$, 
$\mathbb{B}_r \Lambda^{k}(T)\subseteq (\mathbb{P}_r \Lambda^{k}(T)\cap \ker(\tr_{\partial T})).$
Set
\[
\mathcal Q:=\Oplus_{\ell=k}^{d-1}\Oplus_{f\in\Delta_{\ell}(T)}\mathbb B_r\Lambda^k(f).
\]
By \eqref{eq:bubbledec},
$\mathbb P_r\Lambda^k(T)=\mathcal Q\oplus\mathbb B_r\Lambda^k(T)$.
Since $\mathbb B_r\Lambda^k(T)\subseteq\ker(\tr_{\partial T})$, the restricted map
\[
\tr_{\partial T}:\mathcal Q\longrightarrow
\tr_{\partial T}\mathbb P_r\Lambda^k(T)
\]
is surjective.
To prove injectivity, let $\omega\in\mathcal Q$ satisfy $\tr_{\partial T}\omega=0$. For every proper subsimplex $f\subsetneq T$, choose a facet $F\in\Delta_{d-1}(T)$ containing $f$. Transitivity of the trace gives
$\tr_f\omega=\tr_f(\tr_F\omega)=0$.
Thus all DoFs in \eqref{eq:tndofs} assigned to proper faces $f\subsetneq T$ vanish. For fixed $e\in\Delta_s(T)$, this restriction selects exactly the normal degrees below $d-s$ in either normal basis. Hence $\mathcal Q$ is the corresponding partial sum in \eqref{eq:geodec1}. The restricted matrix \eqref{eq:lowertriangular} remains block lower triangular in $\dim e$, with invertible diagonal blocks by the proof of Theorem~\ref{thm:geodec}. Therefore $\omega=0$, proving injectivity and the kernel identity.
\end{proof}

%\mnote{simplify the DoF.}
%Recall that the bubble polynomial satisfies $b_f|_{\partial f} = 0$. In general, for a subsimplex $e\in \Delta(T)$, if $f\cap e^*\neq \varnothing$, then $b_f|_e = 0$. 

%Employing the tensor product, we obtain a straightforward decomposition for the tensor-valued Lagrange element:
%\begin{align}
%\notag \mathbb P_r\Lambda^{k}(T) &:= \mathbb P_r(T)\otimes \alt^k(\mathscr T^T) \\
%\label{intro:polynomialdecomp} &=  \Oplus_{\ell=0}^d\Oplus_{f\in\Delta_{\ell}(T)}\left (\mathbb P_{r-(\ell +1)}(f)b_f\otimes \alt^k(\mathscr T^T)\right ),
%\end{align}
%where
%\begin{itemize}
%\item  $b_f\in \mathbb P_{\ell +1}(f)$ denotes the scalar bubble polynomial on $f$;
%\item $\mathscr T^T$ is the vector space tangent to $T$;
%\item $\alt^k V$ represents the space of alternating multilinear mappings on $V^k$.
%\end{itemize}

\subsection{Finite element spaces}
Let $\Omega\subset\mathbb R^d$ be a bounded Lipschitz polyhedral domain. Let $\{\mathcal T_h\}$ be a family of conforming simplicial partitions of $\Omega$, with $h_T:=\mbox{diam}(T)$ and $h:=\max_{T\in\mathcal T_h}h_T$. The intersection of any two simplices in $\mathcal T_h$ is either empty or a common lower-dimensional subsimplex. Let $\Delta_\ell(\mathcal T_h)$ denote the set of all $\ell$-dimensional subsimplices in $\mathcal T_h$ for $\ell=0,1,\ldots,d$.

%\RV{For every global subsimplex $e$, fix an ordered orientation and a
%specific oriented orthonormal frame of $\mathscr T^e$, independent of
%the element that contains $e$.  For example, one may apply a fixed
%Gram--Schmidt rule to the globally ordered edge vectors.  Orientation
%alone does not choose a frame when $\dim e>1$.  With this extra choice,
%$\dd t_{\sigma}^e\wedge
%\dd_f\widehat{\lambda}_{[f\setminus e]}\in
%\alt^k(\mathscr T^f)$ depends only on the pair $(e,f)$.  We require the
%corresponding DoFs in~\eqref{eq:tndofs} to be single-valued on $f$.}
%That is we look at the $t$-$n$ basis row-wise in Fig.~\ref{fig:tngrid}.

To define the global finite element space $\mathbb P_r\Lambda^k(\mathcal T_h)$, with $r\geq1$, fix an orientation and an oriented orthonormal tangential frame for each global subsimplex $f\in\Delta(\mathcal T_h)$, independent of the element containing $f$. Then
$$
\dd t_\sigma^e\wedge\dd_f\widehat{\lambda}_{[f\setminus e]}\in\alt^k(\mathscr T^f)
$$
depends only on the pair $(e,f)$ with $e\subseteq f$. We assign the corresponding DoF in~\eqref{eq:tndofs} to $f$ and require it to be single-valued across all elements containing $f$. Equivalently, the $t$-$n$ decomposition in Fig.~\ref{fig:tngrid} is read row-wise.%To show the DoFs impose the required continuity to be $H\Lambda^k$-conforming, we c
%an redistribute DoFs \eqref{eq:tndofs} or \eqref{eq:tndofslargek} from $e$ to $f$. Take DoFs \eqref{eq:tndofs} as an example. We can regroup by $f$ not by $e$ and get the following equivalent DoFs
%\RV{Also need to change $\dd  n_{f^*}^f$ to a global dual orientation. For example, for $k=d-1$, $n_F$ is global not local one $\bs n_{f^*}^f$. Will change to the intrinsic one $\dflat \boldsymbol{n}_{e^*\cap f}^{e}$.} \RV{Let $\bs n_{F_i}$ be fixed normal, not the outwards normal direction of $\partial T$.}
%\begin{corollary}\label{thm:redistdofs}
%Let $0\leq k\leq d-1$, and $r\geq1$. On a $d$-dimensional simplex $T$, space $\mathbb P_r\Lambda^{k}(T)$ is uniquely determined by the DoFs: 
%
\begin{theorem}
Let $r\geq1$ and $0\leq k\leq d$, with the global orientations and orthonormal tangential frames fixed above. For each $f\in\Delta_\ell(\mathcal T_h)$, $k\leq\ell\leq d$, associate with $f$ the moments
\begin{equation}
\label{eq:redistdofs1}
\begin{aligned}
\int_e\langle \tr_f(\omega|_T), \dd t_{\sigma}^e\wedge \dd_f \widehat{\lambda}_{[f\setminus  e]} \rangle \ q\dd s,\quad
%& f\in \Delta_{\ell}(\mathcal{T}_h), k\leq \ell \leq d,\\
& e\in\Delta_{s}(f), \ell-k\leq s\leq\ell, \\
& q\in\mathbb P_{r-(s +1)}(e), \sigma\in\Sigma(s+k-\ell,s).
\end{aligned}
\end{equation}
Here $T\in\mathcal T_h$ is any element containing $f$.
Then 
\begin{align*}
\mathbb P_{r}\Lambda^k(\mathcal T_h)&:=\{\omega\in L^2\Lambda^k(\Omega): \omega|_T\in \mathbb P_{r}\Lambda^k(T) \textrm{ for } T\in\mathcal T_h, \\
&\!\!\qquad\textrm{ DoFs \eqref{eq:redistdofs1} are single-valued across } f\in \Delta_{\ell}(\mathcal T_h), \ell = k,\ldots, d-1\}
\end{align*}
coincides with the standard conforming piecewise-polynomial space:
\[
\mathbb P_r\Lambda^k(\mathcal T_h)
=\{\omega\in H\Lambda^k(\Omega):
\omega|_T\in\mathbb P_r\Lambda^k(T)\ \text{for every }T\in\mathcal T_h\}.
\]
\end{theorem}

\begin{proof}
Local unisolvence on each $T$ follows from Theorem~\ref{thm:geodec} by changing the order of the indices $(\ell,s)$ from column-wise to row-wise.

Since $\dd t_\sigma^e\wedge\dd_f\widehat{\lambda}_{[f\setminus e]}\in\alt^k(\mathscr T^f)$ is independent of the element containing $f$, each DoF in~\eqref{eq:redistdofs1} is well-defined on $f$.

If $k=d$, the asserted equality is immediate because $H\Lambda^d(\Omega)=L^2\Lambda^d(\Omega)$ and there are no interelement DoFs. Hence assume $0\leq k\leq d-1$.

Let $\omega\in\mathbb P_r\Lambda^k(\mathcal T_h)$. Consider an interior facet $F=T^+\cap T^-\in\Delta_{d-1}(\mathcal T_h)$. Theorem~\ref{thm:geodec}, applied to $\mathbb P_r\Lambda^k(F)$, shows that each trace $\tr_F(\omega|_{T^\pm})\in\mathbb P_r\Lambda^k(F)$ is uniquely determined by the DoFs associated with all subsimplices $f\subseteq F$. These DoFs are single-valued, so $\tr_F(\omega|_{T^+})=\tr_F(\omega|_{T^-})$. Thus, the tangential trace is single-valued across every interior facet, and hence $\omega\in H\Lambda^k(\Omega)$.

Conversely, let $\omega\in H\Lambda^k(\Omega)$ with $\omega|_T\in\mathbb P_r\Lambda^k(T)$ for every $T\in\mathcal T_h$. Its tangential trace is single-valued across every interior facet. Since traces are pullbacks and hence transitive under restriction to lower-dimensional subsimplices, these facet traces induce the same trace on each $f\in\Delta_\ell(\mathcal T_h)$, $k\leq\ell\leq d-1$, from every element containing $f$; equivalently, one may propagate the equality through the facet-connected star of $f$. Therefore every DoF in \eqref{eq:redistdofs1} is single-valued, so $\omega\in\mathbb P_r\Lambda^k(\mathcal T_h)$.
\end{proof}

\begin{remark}\rm 
An equivalent description of the finite element space uses the following moments associated with each $f\in\Delta_\ell(\mathcal T_h)$, $k\leq\ell\leq\min\{d,r+k-1\}$:
$$
\int_f \tr_f(\omega|_T)\wedge\eta,\qquad \eta\in\mathbb P_{r+k-\ell}^-\Lambda^{\ell-k}(f),
$$
where $T$ contains $f$ and $\mathbb P^-\Lambda$ denotes the first-family polynomial differential forms; see~\cite[(5.1)]{ArnoldFalkWinthe2006Finite}. Relation to \eqref{eq:redistdofs1} through the geometric decomposition of the first family is beyond the scope of this paper.
\end{remark}

\begin{example}\label{ex:bases}
\rm
Using the $t$--$n$ bases in Tables~\ref{fig:tnbasesd3k1} and~\ref{fig:tnbasesd3k2}, we give explicit degrees of freedom from~\eqref{eq:redistdofs1} for $d=3$, with $k=1$ and $k=2$. These are the three-dimensional edge and face elements, respectively.

For $d=3$ and $k=1$, let $f\in\Delta_\ell(\mathcal T_h)$ with $\ell=1,2,3$.
\begin{itemize}[leftmargin=10pt]
\item \textbf{$\ell=1$}. Let $f\in\Delta_1(\mathcal T_h)$ and $e\in\Delta_s(f)$ with $s=0,1$. The associated DoFs are
\begin{align*}
(\boldsymbol v\cdot\boldsymbol t)|_f(\texttt{v}),\quad & e= \texttt{v}\in\Delta_0(f),\\
\int_e\boldsymbol v\cdot\boldsymbol t\,q\,\dd s,\quad &q\in\mathbb P_{r-2}(e),\quad e=f\in\Delta_1(f).
\end{align*}
These DoFs are single-valued because the chosen tangential vector $\boldsymbol t$ is globally defined on $f$. Together, they can be written as
$$
\int_f\boldsymbol v\cdot\boldsymbol t\,q\,\dd s,\qquad q\in\mathbb P_r(f).
$$

\item \textbf{$\ell=2$}. Let $f=F$ and $e\in\Delta_s(F)$ with $s=1,2$. The DoFs are
\begin{align*}
\int_e\boldsymbol v\cdot\nabla_F\lambda_{[F\setminus e]}\,q\,\dd s,\quad &q\in\mathbb P_{r-2}(e),\quad e\in\Delta_1(F),\\
\int_F\boldsymbol v\cdot\boldsymbol t_i^F\,q\,\dd s,\quad &q\in\mathbb P_{r-3}(F),\, i=1,2, \quad e=F\in \Delta_2(F).
\end{align*}
These DoFs are single-valued on each $F\in\Delta_2(\mathcal T_h)$ because $\nabla_F\lambda_{[F\setminus e]}$ and $\boldsymbol t_i^F$ are defined intrinsically on $F$.

\item \textbf{$\ell=3$}. Let $f=T$ and $e\in\Delta_s(T)$ with $s=2,3$. The associated DoFs are
\begin{align*}
\int_F\boldsymbol v\cdot\nabla\lambda_{F^*}\,q\,\dd s,\quad &q\in\mathbb P_{r-3}(F),\quad e=F\in\Delta_2(T),\\
\int_T\boldsymbol v\cdot\boldsymbol q\,\dd x,\quad &\boldsymbol q\in\mathbb P_{r-4}(T;\mathbb R^3), \quad e = T\in \Delta_3(T).
\end{align*}
These DoFs depend on $T$, rather than only on $F$. Thus, the normal component of $\boldsymbol v$ is generally multi-valued across faces. 
%In particular, if two tetrahedra share a face $F$, the corresponding DoFs depend on the chosen tetrahedron.
\end{itemize}

The above DoFs show that the tangential trace on each $f\in\Delta_\ell(\mathcal T_h)$, for $k\leq\ell\leq d-1$, is single-valued. Consequently, the resulting finite element space is a subspace of $H\Lambda^1(\Omega)$.

For $d=3$ and $k=2$, consider $f\in \Delta_{\ell}(\mathcal T_h)$ for $\ell = 2, 3$.
\begin{itemize}[leftmargin =10pt]
\item $\ell = 2$. Then $f = F$ and $e\in \Delta_s(F)$ for $s=0, 1, 2$. The $t$-$n$ basis is given in the first row of Table~\ref{fig:tnbasesd3k2}, and the corresponding DoFs are
\begin{align*}
(\boldsymbol{v}\cdot\boldsymbol{n}_F)|_F(\texttt{v}),\quad & e= \texttt{v}\in\Delta_0(F), \\
\int_e\boldsymbol{v}\cdot\boldsymbol{n}_F \ q\dd s,\quad & q\in\mathbb P_{r-2}(e), \quad e\in\Delta_{1}(F), \\
\int_F\boldsymbol{v}\cdot\boldsymbol{n}_F \ q\dd s,\quad & q\in\mathbb P_{r-3}(F), \quad e = F\in \Delta_{2}(F).
\end{align*}
They are single-valued for each $F\in \Delta_2(\mathcal T_h)$ as $\boldsymbol{n}_F$ depends on $F$ only. By the geometric decomposition of the Lagrange element, these DoFs can be combined into the following set of moments:
$$
\int_F\boldsymbol{v}\cdot\boldsymbol{n}_F \ q\dd s,\quad  q\in\mathbb P_{r}(F). 
$$
\item $\ell = 3$. Then $f = T$ and $e\in \Delta_s(T)$ for $s=1, 2, 3$. The $t$-$n$ basis is given in the second row of Table~\ref{fig:tnbasesd3k2}, and the corresponding DoFs are
\begin{align*}
\int_e\boldsymbol{v}\cdot\boldsymbol{t} \ q\dd s,\quad & q\in\mathbb P_{r-2}(e), \quad e\in\Delta_1(T), \\
\int_F\boldsymbol{v}\cdot\boldsymbol{t}_i^F \ q\dd s,\quad & q\in\mathbb P_{r-3}(F),i=1,2, \quad e=F\in\Delta_2(T), \\
\int_T\boldsymbol{v}\cdot\boldsymbol{q}\dd x,\quad & \boldsymbol{q}\in\mathbb P_{r-4}(T;\mathbb R^3), \quad e = T\in \Delta_3(T).
\end{align*}
They are single-valued for each $T\in \Delta_{3}(\mathcal T_h)$, which implies the tangential component of $\bs v$ is indeed multi-valued; if two tetrahedra share the same face $F$, DoFs depend on $T$. 
%By the bubble decomposition in Section~\ref{sec:geodecBubbleform},
%these DoFs can also be merged into a single volume moment
%\[
%\int_T \boldsymbol{v} \cdot \boldsymbol{q} \,\dd x,
%\qquad
%\forall\, \boldsymbol{q} \in \mathbb B_r \Lambda^{d-1}(T).
%\]
\end{itemize}
From the above DoFs, we see that the normal trace on each $F\in \Delta_{d-1}(\mathcal T_h)$ is single-valued, and thus the finite element space is a subspace of $H\Lambda^{d-1}(\Omega)$. 
\end{example}

\subsection{de Rham complexes}
We have the following de Rham complex on bubble polynomial forms. For similar results, see~\cite[(2.5)]{Falk2014} and~\cite[Lemma~4.24]{ArnoldFalkWinthe2006Finite}. 
% As the bubble polynomial forms are different, we provide a proof here.
\begin{proposition}\label{lm:bubblecomplex}
For a $d$-dimensional simplex $T$ and $r+k\geq d$, the bubble complex
\begin{equation}\label{eq:bubblesequence}
0 \stackrel{}{\longrightarrow} \mathbb B_{r+k}\Lambda^0 (T) \cdots %\mathbb B_{r+1}\Lambda^{k-1}(T)
\stackrel{\dd}{\longrightarrow} \mathbb B_r\Lambda^k(T)\stackrel{\dd}{\longrightarrow} \mathbb B_{r-1}\Lambda^{k+1}(T)\cdots \stackrel{\dd}{\longrightarrow}\mathbb B_{r+k-d}\Lambda^{d}(T)\stackrel{\int}{\longrightarrow} \mathbb R
\end{equation}
is exact.
\end{proposition}

Although our DoFs in \eqref{eq:redistdofs1} for $\mathbb P_r\Lambda^k(\mathcal T_h)$ differ from those in \cite[(5.1)]{ArnoldFalkWinthe2006Finite}, the preceding theorem shows that they define the same standard conforming piecewise-polynomial space. Therefore the standard finite element de Rham exactness result applies, and we obtain the following complex.

\begin{theorem}[cf. (5.13) in \cite{ArnoldFalkWinthe2006Finite}]
Assume $\Omega$ is contractible. Let $r\geq0$, $0\leq k\leq d$, and $r+k\geq d$. When the terminal polynomial degree is zero, set $\mathbb P_0\Lambda^d(\mathcal T_h):=\{\omega\in L^2\Lambda^d(\Omega):\omega|_T\in\mathbb P_0\Lambda^d(T)\text{ for every }T\in\mathcal T_h\}$.
The finite element de Rham complex
$$
\mathbb R \stackrel{}{\hookrightarrow} \mathbb P_{r+k}\Lambda^0 (\mathcal T_h) \cdots %\mathbb B_{r+1}\Lambda^{k-1}(T)
\stackrel{\dd}{\longrightarrow} \mathbb P_r\Lambda^k(\mathcal T_h)\stackrel{\dd}{\longrightarrow} \mathbb P_{r-1}\Lambda^{k+1}(\mathcal T_h)\cdots \stackrel{\dd}{\longrightarrow}\mathbb P_{r+k-d}\Lambda^{d}(\mathcal T_h)\stackrel{}{\rightarrow} 0
$$
is exact.
\end{theorem}

\subsection{Lagrange-type bases}\label{sec:dualbasis}
Throughout this subsection, let $r\geq1$.
The DoFs in \eqref{eq:tndofs} and \eqref{eq:tndofslargek} are formulated using integrals on subsimplices. The shape functions associated with lower-dimensional subsimplices in the geometric decomposition \eqref{eq:geodec1} may not be dual to the DoFs associated with higher-dimensional subsimplices. In particular, the matrix in \eqref{eq:lowertriangular} is not diagonal. We use a nodal basis of Lagrange type to simplify the construction. 

Recall the simplicial lattice \[
\mathbb{T}_r^d = \left\{ \alpha = (\alpha_0, \alpha_1, \ldots, \alpha_d) \in \mathbb{N}^{0:d} \mid \alpha_0 + \alpha_1 + \ldots + \alpha_d = r \right\}.
\]
A geometric embedding of the algebraic set $\mathbb T^d_r$ is given by
%Define a set of points in $T$ 
\[
\mathcal X_{T} = \left \{\boldsymbol
x_{\alpha} = \frac{1}{r}\sum_{i = 0}^d \alpha_i \texttt{v}_i: \alpha\in \mathbb T^d_r \right \},
\]
where $\{\texttt{v}_i, i=0,\ldots, d\}$ are the vertices of $T$. 
We call $\mathcal X_{T}$ the set of interpolation points, which is known as the principal lattice~\cite{nicolaides1972class}. We refer to~\cite{Chen;Huang:2022FEMcomplex3D} for figures of $\mathcal X_T$ in two and three dimensions. 
For a vertex $\texttt{v}\in\Delta_0(\mathcal T_h)$, set $\mathcal X_{\mathring{\texttt{v}}}:=\{\texttt{v}\}$. The set of interpolation points on a conforming triangulation $\mathcal T_h$ then admits the disjoint union decomposition
%\begin{equation}\label{eq:Xunion}
%\mathcal X_h = \bigcup_{T\in \mathcal T_h}\mathcal X_{T}.
%\end{equation}
%Note that a lot of duplications exist in \eqref{eq:Xunion}. A direct sum is given by 
\begin{equation*}%\label{eq:Xsum}
\mathcal X_h  = \bigsqcup_{\ell = 0}^d \bigsqcup_{f\in \Delta_{\ell}(\mathcal T_h)} \mathcal X_{\mathring{f}},
\end{equation*}
where $\mathcal X_{\mathring{f}}$ denotes the set of interpolation points in the interior of $f$ for $1\leq \dim f\leq d$.
For $0\leq \ell\leq d$, define
\[
\mathring{\mathbb T}_r^\ell
:=
\left\{\beta\in\mathbb N^{\ell+1}: |\beta|=r,\ \beta_i>0,\ i=0,\ldots,\ell\right\}
=
\mathbf 1+\mathbb T_{r-(\ell+1)}^\ell,
\]
where $\mathbf 1=(1,\ldots,1)\in\mathbb N^{\ell+1}$ and the last set is understood to be empty when $r<\ell+1$. Fix a global numbering $\iota(f)$ of all subsimplices $f\in\Delta(\mathcal T_h)$, compatible with their dimensions. Introduce
\begin{align*}
\mathbb T^{d}_r(\mathcal T_h):=\{&\alpha=(\alpha_0,\alpha_1,\ldots,\alpha_d,\alpha_{d+1})\in\mathbb N^{d+2}:  \\
&\text{for some }f\in\Delta_\ell(\mathcal T_h),\quad
(\alpha_0,\ldots,\alpha_\ell)\in\mathring{\mathbb T}_r^\ell,\\
&\alpha_i=0\ \text{for }i=\ell+1,\ldots,d,\quad
\alpha_{d+1}=\iota(f)\}.
\end{align*}
The last component $\alpha_{d+1}$ labels the subsimplex containing the lattice point in its interior. For $\alpha\in\mathbb T_r^d(\mathcal T_h)$, let $f_\alpha\in\Delta_\ell(\mathcal T_h)$ be determined by $\iota(f_\alpha)=\alpha_{d+1}$. For each $T\in\mathcal T_h$ containing $f_\alpha$, let $\alpha^T=(\alpha_0^T,\ldots,\alpha_d^T)\in\mathbb T_r^d$ be the local multi-index obtained by assigning the entries $\alpha_0,\ldots,\alpha_\ell$ to the corresponding vertices of $f_\alpha\subseteq T$ and setting all remaining entries to zero. The corresponding interpolation point is
\[
\boldsymbol x_\alpha=\frac{1}{r}\sum_{i=0}^d\alpha_i^T\texttt{v}_i^T.
\]
This point is independent of the containing element $T$, and the map $\alpha\mapsto\boldsymbol x_\alpha$ is a bijection between $\mathbb T^{d}_r(\mathcal T_h)$ and $\mathcal X_h$.

A basis for the $r$th-order Lagrange finite element space
$V_r^{\rm L}(\mathcal T_h)$ \cite{nicolaides1973class} is given by 
\[
\{\phi_{\alpha} \mid \alpha \in \mathbb T^{d}_r(\mathcal T_h)\}.
\]
The global Lagrange basis function associated with $\boldsymbol x_\alpha$ is defined elementwise by
\[
\phi_\alpha|_T(\boldsymbol x)
=
\begin{cases}
\displaystyle
\frac{1}{\alpha^T!}
\prod_{i=0}^d
\prod_{j=0}^{\alpha_i^T-1}
\bigl(r\lambda_i^T(\boldsymbol x)-j\bigr),
& f_\alpha\subseteq T,\\[3mm]
0,
& f_\alpha\not\subseteq T,
\end{cases}
\]
where $\alpha^T!:=\prod_{i=0}^d\alpha_i^T!$.
This basis is dual to the DoFs defined by nodal interpolation:
\[
\phi_{\alpha}(\boldsymbol x_{\beta})
=
\delta_{\alpha,\beta},
\qquad
\alpha, \beta \in \mathbb T^{d}_r(\mathcal T_h).
\]
We use the globally fixed orientations and oriented orthonormal tangential frames for the pairs $(e,f)$ in the $t$-$n$ bases.
\begin{theorem}
Let $r\geq1$ and $0\leq k\leq d$. For each subsimplex in $\mathcal T_h$, fix an orientation and an oriented orthonormal tangential frame independent of the containing element.
Define the index set
\[
\begin{aligned}
\mathcal I_r^k(\mathcal T_h):=\{(\alpha,e,\sigma,f)\mid &\quad f\in\Delta_\ell(\mathcal T_h),\ e\in\Delta_s(f),\\
&\quad \bs x_\alpha\in\mathcal X_{\mathring e},\ \sigma\in\Sigma(s+k-\ell,s),\\
&\quad k\leq\ell\leq d,\ \ell-k\leq s\leq\ell\}.
\end{aligned}
\]
For $(\alpha,e,\sigma,f)\in\mathcal I_r^k(\mathcal T_h)$, define
\[
N^{\sigma,f}_{\alpha,e}(\omega):=
\left\langle\omega(\bs x_\alpha),\dd t_\sigma^e\wedge\dd_f\widehat\lambda_{[f\setminus e]}\right\rangle.
\]
These DoFs are single-valued across $f$ and uniquely determine $\omega\in\mathbb P_r\Lambda^k(\mathcal T_h)$. A scaled dual basis is $\{\bs\phi_{\alpha,e}^{\sigma,f}\}_{(\alpha,e,\sigma,f)\in\mathcal I_r^k(\mathcal T_h)}$, defined elementwise, for every $T\in\mathcal T_h$, by
\[
\left.\bs\phi_{\alpha,e}^{\sigma,f}\right|_T(\bs x)=
\begin{cases}
\left.\phi_\alpha\right|_T(\bs x)\,\dd t_\sigma^e\wedge\dd\lambda^T_{[f\setminus e]},& f\subseteq T,\\
0,& f\not\subseteq T,
\end{cases}
\]
where $\lambda_i^T$ are the barycentric coordinates on $T$.
\end{theorem}
\begin{proof}
On each element, scalar nodal duality and Corollary~\ref{cor:dualtnbasis} give a diagonal DoF--basis pairing with nonzero diagonal entry $c_{e,f}:=\prod_{i\in f\setminus e}|\nabla_{e\cup\{i\}}\lambda_i|^2$, where an empty product equals one. This factor depends only on $(e,f)$.

To check conformity, consider a facet $F$ of an element containing $f$. If $e\not\subseteq F$, then $\phi_\alpha|_F=0$. If $e\subseteq F$ but $f\not\subseteq F$, some $i\in f\setminus F$ contributes a factor $\dd\lambda_i^T$ with zero trace on $F$. If $f\subseteq F$, all factors have the same trace from either adjacent element. The functions therefore assemble conformingly with support on the elements containing $f$. Their global pairing is diagonal, and local unisolvence proves the asserted basis property. Dividing each basis function by $c_{e,f}$ gives an exactly dual basis.
\end{proof}

In~\cite{ChenChenHuangWei2024}, we have implemented high-order edge and face elements in two and three dimensions using the Lagrange basis. For example, for $d = 3$ and $k = 1$, the $t$-$n$ basis for degrees of freedom is given in Table~\ref{fig:tnbasesd3k1}, and the shape function basis is based on the $t$-$n$ basis in Table~\ref{fig:tnbasesd3k2}. Coupled with the standard Lagrange basis and nodal evaluation at the lattice points, we obtain dual bases for the second family of N\'{e}d\'{e}lec edge elements~\cite{Nedelec1986}. DoFs are single-valued on each face $f$ with $\dim f = \ell$ to impose the desired tangential continuity; see Example \ref{ex:bases}.

\section{Conclusion}\label{sec:conclusion}
Building on novel tangential-normal ($t$-$n$) decompositions adapted to subsimplices, we give new algebraic structures and explicit bases for the second family of finite element differential forms that address key computational challenges in FEEC. The geometric decompositions we have constructed provide a systematic, dimension-independent approach for constructing local bases on simplicial meshes. 

The Lagrange-type bases developed here are built on the simplicial lattice associated with $\mathbb P_r(T)$. In future work, we will extend this construction to the first family of finite element differential forms using the interior lattice points associated with $\mathbb P_{r+k+1}(T)$.

\bibliographystyle{abbrv}
\bibliography{FEEC}
\end{document}